\documentclass[letterpaper,english]{article}
\usepackage[T1]{fontenc}
\usepackage[latin9]{inputenc}
\usepackage{longtable}
\usepackage{amsmath}

\makeatletter

\pdfpageheight\paperheight
\pdfpagewidth\paperwidth

\providecommand{\LyX}{L\kern-.1667em\lower.25em\hbox{Y}\kern-.125emX\@}
\providecommand{\tabularnewline}{\\}

\makeatother

\usepackage{babel}
\begin{document}

\title{Induction and Analogy\\
in a Problem of Finite Sums\thanks{This work is licensed under the CC BY 4.0, a Creative Commons Attribution
License.}}

\author{Ryan Zielinski\\
ryan\_zielinski@fastmail.com}

\date{1 July 2016\thanks{Updated on 23 July 2016 to correct typos, turn suggestions into formal
exercises, and clarify the discussion in the last section of the Appendix.} \thanks{Updated on 9 January 2019. See the ``Note to the Reader.''}}
\maketitle
\begin{description}
\item [{Abstract}] What is a general expression for 
\[
\sum_{k=1}^{n}k^{m}=1^{m}+2^{m}+3^{m}+\cdots+n^{m},
\]
where $m$ is a positive integer? Answering this question will be
the aim of the paper. We assume the reader has some mathematical maturity
and is familiar with the results for 
\[
1+2+3+\cdots+n,
\]
\[
1^{2}+2^{2}+3^{2}+\cdots+n^{2},
\]
 and 
\[
1^{3}+2^{3}+3^{3}+\cdots+n^{3}
\]
 and wants to discover the expressions for higher powers. However,
we will take the unorthodox approach of presenting the material from
the point of view of someone who is trying to solve the problem himself.
\item [{Keywords}] analogy, Johann Faulhaber, finite sums, heuristics,
inductive reasoning, number theory, George Polya, problem solving,
teaching of mathematics
\end{description}
\pagebreak{}

(left blank intentionally)\pagebreak{}

\section*{Note to the Reader}

By writing the paper ``Faulhaber and Bernoulli'' the author believes
he has fulfilled the obligation imposed by Section 23, ``Final Remarks,''
of ``Induction and Analogy in a Problem of Finite Sums.'' Therefore
he has decided to revise the earlier work.

Most of the changes have been made to the appendix. The final section
and all references to it have been removed; ``Faulhaber and Bernoulli''
now supercedes such material. In the second section an expression
has been added. The appendix, as a whole, has been formatted appropriately.

Throughout the text, typos have been corrected. A few stylistic changes,
so inconsequential that the reader might fail to notice them, have
been introduced. Otherwise, hardly a word has been changed. The author
has preserved the original work.\pagebreak{}

\tableofcontents{}

\pagebreak{}

\part*{Part 1}

\section{Introduction}

What is a general expression for 
\[
\sum_{k=1}^{n}k^{m}=1^{m}+2^{m}+3^{m}+\cdots+n^{m},
\]
where $m$ is a positive integer? A passage in the Preface to Volume
1 of \cite{key-2} encapsulates the approach we will take to answer
this question:
\begin{quotation}
Mathematics is regarded as a demonstrative science. Yet, this is only
one of its aspects. Finished mathematics presented in a finished form
appears as purely demonstrative, consisting of proofs only. Yet mathematics
in the making resembles any other human knowledge in the making. You
have to guess a mathematical theorem before you prove it; you have
to guess the idea of the proof before you carry through the details.
You have to combine observations and follow analogies; you have to
try and try again. The result of the mathematician's creative work
is demonstrative reasoning, a proof; but the proof is discovered by
plausible reasoning, by guessing. If the learning of mathematics reflects
to any degree the invention of mathematics, it must have a place for
guessing, for plausible inference.
\end{quotation}
We start with the simplest case and a historical example.

\section{A Proof from Antiquity}

What is the sum of $1+2+3+\cdots+n$? Can we write an expression for
it in terms of $n$? What about the special case of $1+2+3+4+5$?
A proof from antiquity proceeds as follows.\footnote{See the article ``Induction and Mathematical Induction'' of \cite{key-1}.}

\bigskip{}
Suppose we draw the following diagram: %
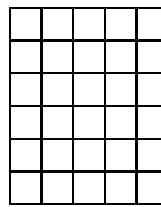
\begin{tabular}{|c|c|c|c|c|}
\hline 
 &  &  &  & \tabularnewline
\hline 
 &  &  &  & \tabularnewline
\hline 
 &  &  &  & \tabularnewline
\hline 
 &  &  &  & \tabularnewline
\hline 
 &  &  &  & \tabularnewline
\hline 
 &  &  &  & \tabularnewline
\hline 
\end{tabular}. There \bigskip{}
are 5 columns and 6 rows . If we think of each square as having unit
area then the total area of the squares is equal to $5\cdot6$. Looking
at the total area in a different way, we may divide it into $1+2+3+4+5$
squares and $1+2+3+4+5$ squares:\footnote{The diagram is not to scale.}\bigskip{}
\begin{tabular}{|c|c|c|c|c|}
\hline 
 &  &  &  & \tabularnewline
\hline 
X &  &  &  & \tabularnewline
\hline 
X & X &  &  & \tabularnewline
\hline 
X & X & X &  & \tabularnewline
\hline 
X & X & X & X & \tabularnewline
\hline 
X & X & X & X & X\tabularnewline
\hline 
\end{tabular}. Therefore we have \bigskip{}
\[
\left(1+2+3+4+5\right)+\left(1+2+3+4+5\right)=5\cdot6,
\]
which we may rewrite as

\[
1+2+3+4+5=\frac{5\cdot6}{2}.
\]
 We generalize that 
\begin{equation}
\sum_{k=1}^{n}k=1+2+3+\cdots+n=\frac{n\left(n+1\right)}{2},\label{eq:1}
\end{equation}
 which is the well-known result.

\bigskip{}
What about higher powers? What is 
\[
\sum_{k=1}^{n}k^{2}=1^{2}+2^{2}+3^{2}+\cdots+n^{2}\,?
\]
 What is 
\[
\sum_{k=1}^{n}k^{3}=1^{3}+2^{3}+3^{3}+\cdots+n^{3}\,?
\]
 Let us postpone looking at these cases and turn to sums of fourth
powers.

\bigskip{}
What is 
\[
\sum_{k=1}^{n}k^{4}=1^{4}+2^{4}+3^{4}+\cdots+n^{4}\,?
\]
 Can we write an expression for it in terms of $n$? Let us look at
the particular case of 
\[
\sum_{k=1}^{5}k^{4}=1^{4}+2^{4}+3^{4}+4^{4}+5^{4}
\]
 and try to adopt the proof from antiquity.

\bigskip{}
Suppose we write %
\begin{tabular}{|c|c|c|c|c|}
\hline 
$1^{3}$ & $2^{3}$ & $3^{3}$ & $4^{3}$ & $5^{3}$\tabularnewline
\hline 
$1^{3}$ & $2^{3}$ & $3^{3}$ & $4^{3}$ & $5^{3}$\tabularnewline
\hline 
$1^{3}$ & $2^{3}$ & $3^{3}$ & $4^{3}$ & $5^{3}$\tabularnewline
\hline 
$1^{3}$ & $2^{3}$ & $3^{3}$ & $4^{3}$ & $5^{3}$\tabularnewline
\hline 
$1^{3}$ & $2^{3}$ & $3^{3}$ & $4^{3}$ & $5^{3}$\tabularnewline
\hline 
$1^{3}$ & $2^{3}$ & $3^{3}$ & $4^{3}$ & $5^{3}$\tabularnewline
\hline 
\end{tabular} . Notice that 
\[
1^{4}+2^{4}+3^{4}+4^{4}+5^{4}=1\cdot1^{3}+2\cdot2^{3}+3\cdot3^{3}+4\cdot4^{3}+5\cdot5^{3}.
\]
\bigskip{}
 On the top part of the table we have written those terms : %
\begin{tabular}{|c|c|c|c|c|}
\hline 
$1^{3}$ & $2^{3}$ & $3^{3}$ & $4^{3}$ & $5^{3}$\tabularnewline
\hline 
 & $2^{3}$ & $3^{3}$ & $4^{3}$ & $5^{3}$\tabularnewline
\hline 
 &  & $3^{3}$ & $4^{3}$ & $5^{3}$\tabularnewline
\hline 
 &  &  & $4^{3}$ & $5^{3}$\tabularnewline
\hline 
 &  &  &  & $5^{3}$\tabularnewline
\hline 
 &  &  &  & \tabularnewline
\hline 
\end{tabular} . On the bottom part of the table we \bigskip{}
have written %
\begin{tabular}{|c|c|c|c|c|}
\hline 
 &  &  &  & \tabularnewline
\hline 
$1^{3}$ &  &  &  & \tabularnewline
\hline 
$1^{3}$ & $2^{3}$ &  &  & \tabularnewline
\hline 
$1^{3}$ & $2^{3}$ & $3^{3}$ &  & \tabularnewline
\hline 
$1^{3}$ & $2^{3}$ & $3^{3}$ & $4^{3}$ & \tabularnewline
\hline 
$1^{3}$ & $2^{3}$ & $3^{3}$ & $4^{3}$ & $5^{3}$\tabularnewline
\hline 
\end{tabular} , which is 
\[
1^{3}+\left(1^{3}+2^{3}\right)+\left(1^{3}+2^{3}+3^{3}\right)+\left(1^{3}+2^{3}+3^{3}+4^{3}\right)+\left(1^{3}+2^{3}+3^{3}+4^{3}+5^{3}\right).
\]
 We may rewrite it as 
\[
\sum_{k=1}^{5}\sum_{l=1}^{k}l^{3}.
\]
 Therefore the sum of the entire table is equal to 
\[
\sum_{k=1}^{5}k^{4}+\sum_{k=1}^{5}\sum_{l=1}^{k}l^{3}.
\]

\bigskip{}
For the next step, remember that in the original proof we wrote the
sum of the area in two different ways. In this proof we can write
the sum of the terms in the diagram in a second way too: 
\[
6\cdot\sum_{k=1}^{5}k^{3}.
\]
 If we set the sums equal to one another then we get 
\[
\sum_{k=1}^{5}k^{4}+\sum_{k=1}^{5}\sum_{l=1}^{k}l^{3}=6\cdot\sum_{k=1}^{5}k^{3},
\]
 which we may rewrite as 
\[
\sum_{k=1}^{5}k^{3+1}=\left(5+1\right)\cdot\sum_{k=1}^{5}k^{3}-\sum_{k=1}^{5}\sum_{l=1}^{k}l^{3}.
\]
 We generalize that 
\begin{equation}
\sum_{k=1}^{n}k^{3+1}=\left(n+1\right)\cdot\sum_{k=1}^{n}k^{3}-\sum_{k=1}^{n}\sum_{l=1}^{k}l^{3}.\label{eq:2}
\end{equation}

\bigskip{}
We have reduced a sum of terms each raised to the fourth power to
a sum of terms each raised to the third power. How do we simplify
such expressions?

\section{The First Catch\label{sub:3 The First Catch}}

Before we try to simplify expression \ref{eq:2} and other ones like
it, let us pause to prove it in the general case. It will be easier
to write it as 
\begin{equation}
\sum_{k=1}^{n}k^{m+1}+\sum_{k=1}^{n}\sum_{l=1}^{k}l^{m}=\left(n+1\right)\cdot\sum_{k=1}^{n}k^{m},\label{eq:3}
\end{equation}
where $m$ is a fixed, positive integer.

\bigskip{}
We proceed by mathematical induction. Previously we established the
result for $\sum_{k=1}^{5}k^{4}$. Let us assume the result is true
for some $n\geq5$ and a fixed, positive integer $m$. Then we may
write 
\[
\sum_{k=1}^{n+1}k^{m+1}+\sum_{k=1}^{n+1}\sum_{l=1}^{k}l^{m}=\sum_{k=1}^{n}k^{m+1}+\left(n+1\right)^{m+1}+\sum_{k=1}^{n}\sum_{l=1}^{k}l^{m}+\sum_{l=1}^{n+1}l^{m}
\]
\[
=\sum_{k=1}^{n}k^{m+1}+\sum_{k=1}^{n}\sum_{l=1}^{k}l^{m}+\left(n+1\right)^{m+1}+\sum_{l=1}^{n+1}l^{m}
\]
\[
=\left(n+1\right)\cdot\sum_{k=1}^{n}k^{m}+\left(n+1\right)^{m+1}+\sum_{l=1}^{n+1}l^{m}
\]
\[
=\left(n+1\right)\cdot\left(\sum_{k=1}^{n}k^{m}+\left(n+1\right)^{m}\right)+\sum_{l=1}^{n+1}l^{m}
\]
\[
=\left(n+1\right)\cdot\sum_{k=1}^{n+1}k^{m}+\sum_{l=1}^{n+1}l^{m}.
\]
 Notice that $\sum_{k=1}^{n+1}k^{m}=\sum_{l=1}^{n+1}l^{m}$. The same
sum is expressed in two different notations. Therefore we may write
\[
\sum_{k=1}^{n+1}k^{m+1}+\sum_{k=1}^{n+1}\sum_{l=1}^{k}l^{m}=\left(n+2\right)\cdot\sum_{k=1}^{n+1}k^{m}.
\]
 We have proved the desired expression.

\section{Back to the Hunt\label{sec:Back-to-the}}

Earlier, we were interested in rewriting expression \ref{eq:2}, 
\[
\sum_{k=1}^{n}k^{3+1}=\left(n+1\right)\cdot\sum_{k=1}^{n}k^{3}-\sum_{k=1}^{n}\sum_{l=1}^{k}l^{3},
\]
 into something simpler, something more along the lines of expression
\ref{eq:1}: 
\[
\sum_{k=1}^{n}k=1+2+3+\cdots+n=\frac{n\left(n+1\right)}{2}.
\]
 Along the way we mentioned the sums of squares and cubes, but did
not try to find expressions for them. Let us do that now.

\bigskip{}
What is 
\[
\sum_{k=1}^{n}k^{2}=1^{2}+2^{2}+3^{2}+\cdots+n^{2}\,?
\]
 Expression \ref{eq:3}, which we proved in Section \ref{sub:3 The First Catch},
tells us 
\begin{equation}
\sum_{k=1}^{n}k^{2}=\left(n+1\right)\cdot\sum_{k=1}^{n}k-\sum_{k=1}^{n}\sum_{l=1}^{k}l.\label{eq:4}
\end{equation}
 We know also that 
\[
\sum_{k=1}^{n}k=\frac{n\left(n+1\right)}{2}=\frac{n+n^{2}}{2}.
\]
 Suppose we \textbf{\textit{assume}} that 
\[
\sum_{k=1}^{n}\sum_{l=1}^{k}l=\frac{\sum_{k=1}^{n}k+\sum_{k=1}^{n}k^{2}}{2}.
\]
 Then we may rewrite expression \ref{eq:4} as 
\[
\sum_{k=1}^{n}k^{2}=\left(n+1\right)\cdot\sum_{k=1}^{n}k-\frac{\sum_{k=1}^{n}k+\sum_{k=1}^{n}k^{2}}{2}
\]
\[
\sum_{k=1}^{n}k^{2}=\left(n+1\right)\cdot\sum_{k=1}^{n}k-\frac{1}{2}\cdot\sum_{k=1}^{n}k-\frac{1}{2}\cdot\sum_{k=1}^{n}k^{2}
\]
\[
\frac{3}{2}\cdot\sum_{k=1}^{n}k^{2}=\left(n+\frac{1}{2}\right)\cdot\sum_{k=1}^{n}k,
\]
 which is 
\[
\sum_{k=1}^{n}k^{2}=\frac{2}{3}\cdot\frac{2n+1}{2}\cdot\sum_{k=1}^{n}k,
\]
 which is the familiar 
\begin{equation}
\sum_{k=1}^{n}k^{2}=\frac{2n+1}{3}\cdot\frac{n\left(n+1\right)}{2}.\label{eq:5}
\end{equation}

\bigskip{}
What is 
\[
\sum_{k=1}^{n}k^{3}=1^{3}+2^{3}+3^{3}+\cdots+n^{3}\,?
\]
 Expression \ref{eq:3} tells us 
\begin{equation}
\sum_{k=1}^{n}k^{3}=\left(n+1\right)\cdot\sum_{k=1}^{n}k^{2}-\sum_{k=1}^{n}\sum_{l=1}^{k}l^{2}.\label{eq:6}
\end{equation}
 We know also that 
\[
\sum_{k=1}^{n}k^{2}=\frac{2n+1}{3}\cdot\frac{n\left(n+1\right)}{2}=\frac{n+3n^{2}+2n^{3}}{6}.
\]
 Suppose we \textbf{\textit{assume}} that 
\[
\sum_{k=1}^{n}\sum_{l=1}^{k}l^{2}=\frac{\sum_{k=1}^{n}k+3\cdot\sum_{k=1}^{n}k^{2}+2\cdot\sum_{k=1}^{n}k^{3}}{6}.
\]
 Then we may rewrite expression \ref{eq:6} as 
\[
\sum_{k=1}^{n}k^{3}=\left(n+1\right)\cdot\sum_{k=1}^{n}k^{2}-\frac{\sum_{k=1}^{n}k+3\cdot\sum_{k=1}^{n}k^{2}+2\cdot\sum_{k=1}^{n}k^{3}}{6}
\]
\[
\sum_{k=1}^{n}k^{3}=\left(n+1\right)\cdot\sum_{k=1}^{n}k^{2}-\frac{1}{6}\cdot\sum_{k=1}^{n}k-\frac{1}{2}\cdot\sum_{k=1}^{n}k^{2}-\frac{1}{3}\cdot\sum_{k=1}^{n}k^{3}
\]
\[
\frac{4}{3}\cdot\sum_{k=1}^{n}k^{3}=\left(n+\frac{1}{2}\right)\cdot\sum_{k=1}^{n}k^{2}-\frac{1}{6}\cdot\sum_{k=1}^{n}k
\]
\[
\sum_{k=1}^{n}k^{3}=\frac{3}{4}\cdot\frac{2n+1}{2}\cdot\frac{2n+1}{3}\cdot\frac{n\left(n+1\right)}{2}-\frac{3}{24}\cdot\frac{n\left(n+1\right)}{2}
\]
\[
=\frac{\left(2n+1\right)^{2}}{8}\cdot\frac{n\left(n+1\right)}{2}-\frac{1}{8}\cdot\frac{n\left(n+1\right)}{2}
\]
\[
=\frac{1}{8}\cdot\left(\left(2n+1\right)^{2}-1\right)\cdot\frac{n\left(n+1\right)}{2},
\]
 which is 
\[
\sum_{k=1}^{n}k^{3}=\frac{1}{8}\cdot\left(4n^{2}+4n\right)\cdot\frac{n\left(n+1\right)}{2},
\]
 which is the familiar 
\begin{equation}
\sum_{k=1}^{n}k^{3}=\left(\frac{n\left(n+1\right)}{2}\right)^{2}.\label{eq:7}
\end{equation}
 What's going on?

\section{An Explanation\label{sub:5 An Explanation}}

Let us look at the assumptions we made in Section \ref{sec:Back-to-the}.
Consider the expression for the sum of squares: 
\[
\sum_{k=1}^{n}k^{2}=\frac{2n+1}{3}\cdot\frac{n\left(n+1\right)}{2}=\frac{n+3n^{2}+2n^{3}}{6}.
\]
 Suppose we rewrite $\frac{n+3n^{2}+2n^{3}}{6}$ as $\frac{1}{6}\cdot n+\frac{3}{6}\cdot n^{2}+\frac{2}{6}\cdot n^{3}$.
In an analogous fashion we may write 
\[
\sum_{l=1}^{k}l^{2}=\frac{1}{6}\cdot k+\frac{3}{6}\cdot k^{2}+\frac{2}{6}\cdot k^{3},
\]
 which implies 
\[
\sum_{k=1}^{n}\sum_{l=1}^{k}l^{2}=\sum_{k=1}^{n}\left(\frac{1}{6}\cdot k+\frac{3}{6}\cdot k^{2}+\frac{2}{6}\cdot k^{3}\right)
\]
\[
=\frac{1}{6}\cdot\sum_{k=1}^{n}k+\frac{3}{6}\cdot\sum_{k=1}^{n}k^{2}+\frac{2}{6}\cdot\sum_{k=1}^{n}k^{3}
\]
\[
=\frac{\sum_{k=1}^{n}k+3\cdot\sum_{k=1}^{n}k^{2}+2\cdot\sum_{k=1}^{n}k^{3}}{6},
\]
 which is analogous to $\sum_{k=1}^{n}k^{2}=\frac{n+3n^{2}+2n^{3}}{6}$
.

\bigskip{}
We have justified the previous, unproven steps. More important, if
we have an expression for $\sum_{l=1}^{k}l^{m}$ then we can derive
an expression for $\sum_{k=1}^{n}\sum_{l=1}^{k}l^{m}$ .

\section{Now for the Big Game\label{sec:Now-for-the}}

Let us put together our previous work to derive expressions we might
not have seen before. What is 
\[
\sum_{k=1}^{n}k^{4}=1^{4}+2^{4}+3^{4}+\cdots+n^{4}\,?
\]
 Expression \ref{eq:3} tells us 
\begin{equation}
\sum_{k=1}^{n}k^{4}=\left(n+1\right)\cdot\sum_{k=1}^{n}k^{3}-\sum_{k=1}^{n}\sum_{l=1}^{k}l^{3}.\label{eq:8}
\end{equation}
 We know also that 
\[
\sum_{k=1}^{n}k^{3}=\left(\frac{n\left(n+1\right)}{2}\right)^{2}=\frac{n^{2}+2n^{3}+n^{4}}{4}.
\]
 The discussion in Section \ref{sub:5 An Explanation} tells us 
\[
\sum_{k=1}^{n}\sum_{l=1}^{k}l^{3}=\frac{\sum_{k=1}^{n}k^{2}+2\cdot\sum_{k=1}^{n}k^{3}+\sum_{k=1}^{n}k^{4}}{4}.
\]
 Therefore we may rewrite expression \ref{eq:8} as 

\[
\sum_{k=1}^{n}k^{4}=\left(n+1\right)\cdot\sum_{k=1}^{n}k^{3}-\frac{\sum_{k=1}^{n}k^{2}+2\cdot\sum_{k=1}^{n}k^{3}+\sum_{k=1}^{n}k^{4}}{4}
\]
\[
\sum_{k=1}^{n}k^{4}=\left(n+1\right)\cdot\sum_{k=1}^{n}k^{3}-\frac{1}{4}\cdot\sum_{k=1}^{n}k^{2}-\frac{1}{2}\cdot\sum_{k=1}^{n}k^{3}-\frac{1}{4}\cdot\sum_{k=1}^{n}k^{4}
\]
\[
\frac{5}{4}\cdot\sum_{k=1}^{n}k^{4}=\frac{2n+1}{2}\cdot\sum_{k=1}^{n}k^{3}-\frac{1}{4}\cdot\sum_{k=1}^{n}k^{2}
\]
\[
\sum_{k=1}^{n}k^{4}=\frac{4}{5}\cdot\frac{2n+1}{2}\cdot\left(\frac{n\left(n+1\right)}{2}\right)^{2}-\frac{4}{5}\cdot\frac{1}{4}\cdot\frac{2n+1}{3}\cdot\frac{n\left(n+1\right)}{2}
\]
\[
=\frac{3}{3}\cdot\frac{4}{5}\cdot\frac{2n+1}{2}\cdot\left(\frac{n\left(n+1\right)}{2}\right)^{2}-\frac{1}{5}\cdot\frac{2n+1}{3}\cdot\frac{n\left(n+1\right)}{2}
\]
\[
=\frac{2n+1}{3}\cdot\frac{12}{10}\cdot\left(\frac{n\left(n+1\right)}{2}\right)^{2}-\frac{1}{5}\cdot\frac{2n+1}{3}\cdot\frac{n\left(n+1\right)}{2},
\]
 which is a less familiar 
\[
\sum_{k=1}^{n}k^{4}=\frac{2n+1}{3}\cdot\left(\frac{6}{5}\cdot\left(\frac{n\left(n+1\right)}{2}\right)^{2}-\frac{1}{5}\cdot\frac{n\left(n+1\right)}{2}\right)
\]
 or 
\[
\sum_{k=1}^{n}k^{4}=\frac{1}{5}\cdot\left(6\cdot\frac{n\left(n+1\right)}{2}-1\right)\cdot\frac{2n+1}{3}\cdot\frac{n\left(n+1\right)}{2}
\]
or 
\begin{equation}
\sum_{k=1}^{n}k^{4}=\frac{3n\left(n+1\right)-1}{5}\cdot\frac{2n+1}{3}\cdot\frac{n\left(n+1\right)}{2},\label{eq:9}
\end{equation}
 which simplifies to 
\[
\sum_{k=1}^{n}k^{4}=\frac{-n+10n^{3}+15n^{4}+6n^{5}}{30}.
\]

\bigskip{}
What is 
\[
\sum_{k=1}^{n}k^{5}=1^{5}+2^{5}+3^{5}+\cdots+n^{5}\,?
\]
 Expression \ref{eq:3} tells us 
\begin{equation}
\sum_{k=1}^{n}k^{5}=\left(n+1\right)\cdot\sum_{k=1}^{n}k^{4}-\sum_{k=1}^{n}\sum_{l=1}^{k}l^{4}.\label{eq:10}
\end{equation}
 We just derived 
\[
\sum_{k=1}^{n}k^{4}=\frac{-n+10n^{3}+15n^{4}+6n^{5}}{30}.
\]
 The discussion in Section \ref{sub:5 An Explanation} tells us 
\[
\sum_{k=1}^{n}\sum_{l=1}^{k}l^{4}=\frac{-\sum_{k=1}^{n}k+10\cdot\sum_{k=1}^{n}k^{3}+15\cdot\sum_{k=1}^{n}k^{4}+6\cdot\sum_{k=1}^{n}k^{5}}{30}.
\]
 Therefore we may rewrite expression \ref{eq:10} as 
\[
\sum_{k=1}^{n}k^{5}=\left(n+1\right)\cdot\sum_{k=1}^{n}k^{4}-\frac{-\sum_{k=1}^{n}k+10\cdot\sum_{k=1}^{n}k^{3}+15\cdot\sum_{k=1}^{n}k^{4}+6\cdot\sum_{k=1}^{n}k^{5}}{30}
\]
\[
\sum_{k=1}^{n}k^{5}=(n+1)\cdot\sum_{k=1}^{n}k^{4}+\frac{1}{30}\cdot\sum_{k=1}^{n}k-\frac{1}{3}\cdot\sum_{k=1}^{n}k^{3}-\frac{1}{2}\cdot\sum_{k=1}^{n}k^{4}-\frac{1}{5}\cdot\sum_{k=1}^{n}k^{5}
\]
\[
\frac{6}{5}\cdot\sum_{k=1}^{n}k^{5}=\frac{2n+1}{2}\cdot\sum_{k=1}^{n}k^{4}+\frac{1}{30}\cdot\sum_{k=1}^{n}k-\frac{1}{3}\cdot\sum_{k=1}^{n}k^{3}.
\]

\bigskip{}
For $\sum_{k=1}^{n}k^{4}$ we substitute 
\[
\frac{2n+1}{3}\cdot\left(\frac{6}{5}\cdot\left(\frac{n\left(n+1\right)}{2}\right)^{2}-\frac{1}{5}\cdot\frac{n\left(n+1\right)}{2}\right)
\]
 and rewrite the expression as 
\begin{eqnarray*}
\sum_{k=1}^{n}k^{5} & = & \frac{5}{6}\cdot\frac{2n+1}{2}\cdot\frac{2n+1}{3}\cdot\left(\frac{6}{5}\cdot\left(\frac{n\left(n+1\right)}{2}\right)^{2}-\frac{1}{5}\cdot\frac{n\left(n+1\right)}{2}\right)\\
 &  & +\frac{5}{6}\cdot\frac{1}{30}\cdot\frac{n\left(n+1\right)}{2}-\frac{5}{6}\cdot\frac{1}{3}\cdot\left(\frac{n\left(n+1\right)}{2}\right)^{2}.
\end{eqnarray*}
 For the case of $\sum_{k=1}^{n}k^{3}$ we rounded up terms of the
form $(2n+1)$. We do the same here: 
\begin{eqnarray*}
\sum_{k=1}^{n}k^{5} & = & \left(2n+1\right)^{2}\cdot\frac{5}{6}\cdot\frac{1}{2}\cdot\frac{1}{3}\cdot\frac{6}{5}\cdot\left(\frac{n\left(n+1\right)}{2}\right)^{2}-\left(2n+1\right)^{2}\cdot\frac{5}{6}\cdot\frac{1}{2}\cdot\frac{1}{3}\cdot\frac{1}{5}\cdot\frac{n\left(n+1\right)}{2}\\
 &  & +\frac{1}{36}\cdot\frac{n\left(n+1\right)}{2}-\frac{5}{18}\cdot\left(\frac{n\left(n+1\right)}{2}\right)^{2}
\end{eqnarray*}
\[
=\frac{\left(2n+1\right)^{2}}{6}\cdot\left(\frac{n\left(n+1\right)}{2}\right)^{2}-\frac{\left(2n+1\right)^{2}}{36}\cdot\frac{n\left(n+1\right)}{2}+\frac{1}{36}\cdot\frac{n\left(n+1\right)}{2}-\frac{5}{18}\cdot\left(\frac{n\left(n+1\right)}{2}\right)^{2}
\]
\[
=\left(\frac{\left(2n+1\right)^{2}}{6}-\frac{5}{18}\right)\cdot\left(\frac{n\left(n+1\right)}{2}\right)^{2}+\left(\frac{1}{36}-\frac{\left(2n+1\right)^{2}}{36}\right)\cdot\frac{n\left(n+1\right)}{2},
\]
 which is 
\[
\sum_{k=1}^{n}k^{5}=\frac{6n\left(n+1\right)-1}{9}\cdot\left(\frac{n\left(n+1\right)}{2}\right)^{2}-\frac{2n\left(n+1\right)}{18}\cdot\frac{n\left(n+1\right)}{2}.
\]
 We may rewrite it further as 
\[
\sum_{k=1}^{n}k^{5}=\frac{6n\left(n+1\right)-1}{9}\cdot\left(\frac{n\left(n+1\right)}{2}\right)^{2}-\frac{2}{9}\cdot\frac{n\left(n+1\right)}{2}\cdot\frac{n\left(n+1\right)}{2}
\]
\[
\sum_{k=1}^{n}k^{5}=\left(\frac{6n\left(n+1\right)-1}{9}-\frac{2}{9}\right)\cdot\left(\frac{n\left(n+1\right)}{2}\right)^{2}
\]
\[
\sum_{k=1}^{n}k^{5}=\frac{1}{3}\cdot\left(4\cdot\frac{n\left(n+1\right)}{2}-1\right)\cdot\left(\frac{n\left(n+1\right)}{2}\right)^{2}
\]
 or 
\begin{equation}
\sum_{k=1}^{n}k^{5}=\frac{2n\left(n+1\right)-1}{3}\cdot\left(\frac{n\left(n+1\right)}{2}\right)^{2},\label{eq:11}
\end{equation}
 which simplifies to 
\[
\sum_{k=1}^{n}k^{5}=\frac{-n^{2}+5n^{4}+6n^{5}+2n^{6}}{12}.
\]

\pagebreak{}

\section{Summary of Part 1\label{sec:Summary}}

Using the well-known result of 
\[
\sum_{k=1}^{n}k=\frac{n\left(n+1\right)}{2}
\]
 and an analogy of the method from a proof from antiquity we discovered
\[
\sum_{k=1}^{n}k^{m+1}+\sum_{k=1}^{n}\sum_{l=1}^{k}l^{m}=\left(n+1\right)\cdot\sum_{k=1}^{n}k^{m},
\]
 which we proved rigorously. We used the new result to derive 
\[
\sum_{k=1}^{n}k^{2}=\frac{2n+1}{3}\cdot\frac{n\left(n+1\right)}{2}
\]
\[
\sum_{k=1}^{n}k^{3}=\left(\frac{n\left(n+1\right)}{2}\right)^{2}
\]
\begin{eqnarray*}
\sum_{k=1}^{n}k^{4} & = & \frac{3n\left(n+1\right)-1}{5}\cdot\frac{2n+1}{3}\cdot\frac{n\left(n+1\right)}{2}\\
 &  & =\frac{1}{5}\cdot\left(6\cdot\frac{n\left(n+1\right)}{2}-1\right)\cdot\frac{2n+1}{3}\cdot\frac{n\left(n+1\right)}{2}
\end{eqnarray*}
\begin{eqnarray*}
\sum_{k=1}^{n}k^{5} & = & \frac{2n\left(n+1\right)-1}{3}\cdot\left(\frac{n\left(n+1\right)}{2}\right)^{2}\\
 &  & =\frac{1}{3}\cdot\left(4\cdot\frac{n\left(n+1\right)}{2}-1\right)\cdot\left(\frac{n\left(n+1\right)}{2}\right)^{2}.
\end{eqnarray*}

The approach can be continued to derive expressions for $\sum_{k=1}^{n}k^{6}$,
$\sum_{k=1}^{n}k^{7}$, and $\sum_{k=1}^{n}k^{223}$. In other words,
given an expression for $\sum_{k=1}^{n}k^{m}$, we can derive expressions
for 
\[
\sum_{k=1}^{n}k^{m+1},\sum_{k=1}^{n}k^{m+2},\ldots.
\]
 What remains to do is, for any $m$, to find a general expression
for $\sum_{k=1}^{n}k^{m}$.

\pagebreak{}

\part*{Part 2}

\section{Emerging Patterns\label{sec:Emerging-Patterns}}

Do we notice any patterns in the expressions we derived? Let us look
at them again: 
\[
\sum_{k=1}^{n}k^{2}=\frac{2n+1}{3}\cdot\frac{n\left(n+1\right)}{2}
\]
\[
\sum_{k=1}^{n}k^{3}=\left(\frac{n\left(n+1\right)}{2}\right)^{2}
\]
\[
\sum_{k=1}^{n}k^{4}=\frac{3n\left(n+1\right)-1}{5}\cdot\frac{2n+1}{3}\cdot\frac{n\left(n+1\right)}{2}
\]
\[
\sum_{k=1}^{n}k^{5}=\frac{2n\left(n+1\right)-1}{3}\cdot\left(\frac{n\left(n+1\right)}{2}\right)^{2}.
\]
 We see that $\frac{3n\left(n+1\right)-1}{5}$ is the coefficient
for $\sum_{k=1}^{n}k^{4}$ and that $\frac{2n\left(n+1\right)-1}{3}$
is the coefficient for $\sum_{k=1}^{n}k^{5}$. We believe that $\sum_{k=1}^{n}k^{6}$
and $\sum_{k=1}^{n}k^{7}$ will have analogous coefficients, and that
if we place all of the expressions together then we will be able to
discern a general pattern.

\bigskip{}
We said ``coefficients.'' Coefficients for what?\footnote{By ``coefficient'' we mean ``leading term,'' not necessarily a
constant factor.} Let us rewrite the expressions as follows: 
\[
\sum_{k=1}^{n}k^{2}=\frac{2n+1}{3}\cdot\frac{n\left(n+1\right)}{2}
\]
\[
\sum_{k=1}^{n}k^{3}=\left(\frac{n\left(n+1\right)}{2}\right)^{2}
\]
\[
\sum_{k=1}^{n}k^{4}=\frac{3n\left(n+1\right)-1}{5}\cdot\sum_{k=1}^{n}k^{2}
\]
\[
\sum_{k=1}^{n}k^{5}=\frac{2n\left(n+1\right)-1}{3}\cdot\sum_{k=1}^{n}k^{3}.
\]
 We see that $\sum_{k=1}^{n}k^{2}$ appears in the expression for
$\sum_{k=1}^{n}k^{4}$ and that $\sum_{k=1}^{n}k^{3}$ appears in
the expression for $\sum_{k=1}^{n}k^{5}$. Is it possible that 
\[
\sum_{k=1}^{n}k^{6}=E_{6}(n)\cdot\sum_{k=1}^{n}k^{2}
\]
 and 
\[
\sum_{k=1}^{n}k^{7}=O_{7}(n)\cdot\sum_{k=1}^{n}k^{3},
\]
 where $E_{6}(n)$ and $O_{7}(n)$ are yet-to-be-determined coefficients
for $\sum_{k=1}^{n}k^{6}$ and $\sum_{k=1}^{n}k^{7}$? Let us find
out.\footnote{$E$ stands for ``even'' and $O$ stands for ``odd.'' We distinguish
a pattern between even and odd powers.}

\subsection{$\sum_{k=1}^{n}k^{6}$}

What is 
\[
\sum_{k=1}^{n}k^{6}=1^{6}+2^{6}+3^{6}+\cdots+n^{6}\,?
\]
 We know that 
\begin{equation}
\sum_{k=1}^{n}k^{6}=\left(n+1\right)\cdot\sum_{k=1}^{n}k^{5}-\sum_{k=1}^{n}\sum_{l=1}^{k}l^{5}.\label{eq:12}
\end{equation}
 We know also that 
\[
\sum_{k=1}^{n}k^{4}=\frac{3n\left(n+1\right)-1}{5}\cdot\sum_{k=1}^{n}k^{2}
\]
 and 
\begin{eqnarray*}
\sum_{k=1}^{n}k^{5} & = & \frac{2n\left(n+1\right)-1}{3}\cdot\sum_{k=1}^{n}k^{3}\\
 &  & =\frac{-n^{2}+5n^{4}+6n^{5}+2n^{6}}{12},
\end{eqnarray*}
 which implies 
\[
\sum_{k=1}^{n}\sum_{l=1}^{k}l^{5}=\frac{-\sum_{k=1}^{n}k^{2}+5\cdot\sum_{k=1}^{n}k^{4}+6\cdot\sum_{k=1}^{n}k^{5}+2\cdot\sum_{k=1}^{n}k^{6}}{12}.
\]
 Therefore we may rewrite expression \ref{eq:12} as 
\[
\sum_{k=1}^{n}k^{6}=\left(n+1\right)\cdot\sum_{k=1}^{n}k^{5}-\frac{-\sum_{k=1}^{n}k^{2}+5\cdot\sum_{k=1}^{n}k^{4}+6\cdot\sum_{k=1}^{n}k^{5}+2\cdot\sum_{k=1}^{n}k^{6}}{12}
\]
\[
=\left(n+1\right)\cdot\sum_{k=1}^{n}k^{5}+\frac{1}{12}\cdot\sum_{k=1}^{n}k^{2}-\frac{5}{12}\cdot\sum_{k=1}^{n}k^{4}-\frac{1}{2}\cdot\sum_{k=1}^{n}k^{5}-\frac{1}{6}\cdot\sum_{k=1}^{n}k^{6}
\]
\[
\frac{7}{6}\cdot\sum_{k=1}^{n}k^{6}=\left(n+\frac{1}{2}\right)\cdot\sum_{k=1}^{n}k^{5}+\frac{1}{12}\cdot\sum_{k=1}^{n}k^{2}-\frac{5}{12}\cdot\sum_{k=1}^{n}k^{4}
\]
\[
\sum_{k=1}^{n}k^{6}=\frac{6}{7}\cdot\frac{2n+1}{2}\cdot\sum_{k=1}^{n}k^{5}+\frac{1}{14}\cdot\sum_{k=1}^{n}k^{2}-\frac{5}{14}\cdot\sum_{k=1}^{n}k^{4}.
\]

\bigskip{}
For $\sum_{k=1}^{n}k^{4}$ we substitute 
\[
\frac{3n\left(n+1\right)-1}{5}\cdot\sum_{k=1}^{n}k^{2}
\]
 and rewrite the expression as 
\[
\sum_{k=1}^{n}k^{6}=\frac{6}{7}\cdot\frac{2n+1}{2}\cdot\sum_{k=1}^{n}k^{5}+\frac{1}{14}\cdot\sum_{k=1}^{n}k^{2}-\frac{5}{14}\cdot\frac{3n\left(n+1\right)-1}{5}\cdot\sum_{k=1}^{n}k^{2}
\]
\[
=\frac{6}{7}\cdot\frac{2n+1}{2}\cdot\sum_{k=1}^{n}k^{5}+\left(\frac{1}{14}-\frac{5}{14}\cdot\frac{3n\left(n+1\right)-1}{5}\right)\cdot\sum_{k=1}^{n}k^{2}
\]
\[
=\frac{6}{7}\cdot\frac{2n+1}{2}\cdot\sum_{k=1}^{n}k^{5}+\left(\frac{-3n\left(n+1\right)+2}{14}\right)\cdot\sum_{k=1}^{n}k^{2}.
\]
 For $\sum_{k=1}^{n}k^{5}$ we substitute 
\[
\frac{2n\left(n+1\right)-1}{3}\cdot\sum_{k=1}^{n}k^{3}
\]
 and rewrite the left side of the expression as 
\[
\frac{6}{7}\cdot\frac{2n+1}{2}\cdot\frac{2n\left(n+1\right)-1}{3}\cdot\sum_{k=1}^{n}k^{3}=\frac{6}{7}\cdot\frac{2n+1}{3}\cdot\frac{2n\left(n+1\right)-1}{2}\cdot\sum_{k=1}^{n}k^{3}
\]
\[
=\frac{6}{7}\cdot\frac{2n+1}{3}\cdot\frac{2n\left(n+1\right)-1}{2}\cdot\frac{n\left(n+1\right)}{2}\cdot\frac{n\left(n+1\right)}{2}
\]
\[
=\frac{6}{7}\cdot\frac{2n\left(n+1\right)-1}{2}\cdot\frac{n\left(n+1\right)}{2}\cdot\frac{2n+1}{3}\cdot\frac{n\left(n+1\right)}{2}
\]
\[
=\frac{6}{7}\cdot\frac{2n\left(n+1\right)-1}{2}\cdot\frac{n\left(n+1\right)}{2}\cdot\sum_{k=1}^{n}k^{2}.
\]
 Together we have 
\[
\sum_{k=1}^{n}k^{6}=\frac{6}{7}\cdot\frac{2n\left(n+1\right)-1}{2}\cdot\frac{n\left(n+1\right)}{2}\cdot\sum_{k=1}^{n}k^{2}+\left(\frac{-3n\left(n+1\right)+2}{14}\right)\cdot\sum_{k=1}^{n}k^{2}
\]
\[
=\frac{3\left(2n\left(n+1\right)-1\right)\cdot n\left(n+1\right)}{14}\cdot\sum_{k=1}^{n}k^{2}+\left(\frac{-3n\left(n+1\right)+2}{14}\right)\cdot\sum_{k=1}^{n}k^{2}
\]
\[
=\left(\frac{3\left(2n\left(n+1\right)-1\right)\cdot n\left(n+1\right)}{14}+\frac{-3n\left(n+1\right)+2}{14}\right)\cdot\sum_{k=1}^{n}k^{2}
\]
\[
=\left(\frac{6\left(n\left(n+1\right)\right)^{2}-3n\left(n+1\right)-3n\left(n+1\right)+2}{14}\right)\cdot\sum_{k=1}^{n}k^{2},
\]
 which is 
\begin{equation}
\sum_{k=1}^{n}k^{6}=\left(\frac{3\left(n\left(n+1\right)\right)^{2}-3n\left(n+1\right)+1}{7}\right)\cdot\sum_{k=1}^{n}k^{2}.\label{eq:13}
\end{equation}

\bigskip{}
Success! We wrote $\sum_{k=1}^{n}k^{6}$ as 
\[
\sum_{k=1}^{n}k^{6}=E_{6}(n)\cdot\sum_{k=1}^{n}k^{2},
\]
 where 
\[
E_{6}(n)=\left(\frac{3\left(n\left(n+1\right)\right)^{2}-3n\left(n+1\right)+1}{7}\right).
\]
 Another way to write expression \ref{eq:13} is 
\begin{equation}
\sum_{k=1}^{n}k^{6}=\frac{1}{7}\cdot\left(12\cdot\left(\frac{n\left(n+1\right)}{2}\right)^{2}-6\cdot\frac{n\left(n+1\right)}{2}+1\right)\cdot\sum_{k=1}^{n}k^{2}.\label{eq:14}
\end{equation}
 Both expressions simplify to 
\[
\sum_{k=1}^{n}k^{6}=\frac{n-7n^{3}+21n^{5}+21n^{6}+6n^{7}}{42}.
\]

\section{Success}

Let us look at our recent success. We had a hunch we were going to
find an expression like 
\[
\sum_{k=1}^{n}k^{6}=E_{6}(n)\cdot\sum_{k=1}^{n}k^{2},
\]
 and we actually were able to find 
\[
\sum_{k=1}^{n}k^{6}=\left(\frac{3\left(n\left(n+1\right)\right)^{2}-3n\left(n+1\right)+1}{7}\right)\cdot\sum_{k=1}^{n}k^{2}.
\]
 What about $\sum_{k=1}^{n}k^{7}$? Will it turn out the same way?
The calculation in Section \ref{sub:seven} of the Appendix tells
us 
\begin{equation}
\sum_{k=1}^{n}k^{7}=O_{7}(n)\cdot\sum_{k=1}^{n}k^{3},\label{eq:15}
\end{equation}
 where 
\[
O_{7}(n)=\left(\frac{3(n(n+1))^{2}-4n(n+1)+2}{6}\right).
\]
 The same pattern has appeared again. Do we have enough information
to figure out the general case? Let us place all of our results together:
\[
\sum_{k=1}^{n}k^{2}=\frac{2n+1}{3}\cdot\sum_{k=1}^{n}k
\]
\[
\sum_{k=1}^{n}k^{3}=\left(\sum_{k=1}^{n}k\right)^{2}
\]
\[
\sum_{k=1}^{n}k^{4}=\frac{3n\left(n+1\right)-1}{5}\cdot\sum_{k=1}^{n}k^{2}
\]
\[
\sum_{k=1}^{n}k^{5}=\frac{2n\left(n+1\right)-1}{3}\cdot\sum_{k=1}^{n}k^{3}
\]
\[
\sum_{k=1}^{n}k^{6}=\left(\frac{3\left(n\left(n+1\right)\right)^{2}-3n\left(n+1\right)+1}{7}\right)\cdot\sum_{k=1}^{n}k^{2}
\]
\[
\sum_{k=1}^{n}k^{7}=\left(\frac{3(n(n+1))^{2}-4n(n+1)+2}{6}\right)\cdot\sum_{k=1}^{n}k^{3}.
\]
 Can we guess general expressions for $E_{2m}(n)$ and $O_{2m+1}(n)$?

\section{Refining Our Approach\label{sec:Refining-Our-Approach}}

A general pattern still remains a mystery. Perhaps it will help to
calculate $\sum_{k=1}^{n}k^{8}$ and $\sum_{k=1}^{n}k^{9}$. Unfortunately,
we have reached the limits of our approach. We need to modify the
way we carry out the calculations.

\bigskip{}
From the work we have done, we know the major stumbling point in the
calculations has been starting with an expression like 
\[
\sum_{k=1}^{n}k^{6}=\left(n+1\right)\cdot\sum_{k=1}^{n}k^{5}-\frac{-\sum_{k=1}^{n}k^{2}+5\cdot\sum_{k=1}^{n}k^{4}+6\cdot\sum_{k=1}^{n}k^{5}+2\cdot\sum_{k=1}^{n}k^{6}}{12},
\]
 for example. Even with looking for a pattern in $\sum_{k=1}^{n}k^{2}$
or $\sum_{k=1}^{n}k^{3}$, unpacking the terms on the right side has
become too unwieldy. The difficulty is, for each term in the numerator
we have been introducing its coefficients from the start and then
trying to simplify them before moving onto the next term. Let us try
something else.

\subsection{$\sum_{k=1}^{n}k^{8}$}

What is 
\[
\sum_{k=1}^{n}k^{8}=1^{8}+2^{8}+3^{8}+\cdots+n^{8}\,?
\]
 Expression \ref{eq:3}, of Section \ref{sub:3 The First Catch},
tells us 
\begin{equation}
\sum_{k=1}^{n}k^{8}=\left(n+1\right)\cdot\sum_{k=1}^{n}k^{7}-\sum_{k=1}^{n}\sum_{l=1}^{k}l^{7}.\label{eq:16}
\end{equation}
 The calculations for $\sum_{k=1}^{n}k^{7}$ in Section \ref{sub:seven}
of the Appendix tell us 
\[
\sum_{k=1}^{n}k^{7}=\frac{2n^{2}-7n^{4}+14n^{6}+12n^{7}+3n^{8}}{24},
\]
 which implies 
\[
\sum_{k=1}^{n}\sum_{l=1}^{k}l^{7}=\frac{2\cdot\sum_{k=1}^{n}k^{2}-7\cdot\sum_{k=1}^{n}k^{4}+14\cdot\sum_{k=1}^{n}k^{6}+12\cdot\sum_{k=1}^{n}k^{7}+3\cdot\sum_{k=1}^{n}k^{8}}{24}.
\]
 Therefore we may rewrite expression \ref{eq:16} as 
\[
\sum_{k=1}^{n}k^{8}=\left(n+1\right)\cdot\sum_{k=1}^{n}k^{7}-\frac{2\cdot\sum_{k=1}^{n}k^{2}-7\cdot\sum_{k=1}^{n}k^{4}+14\cdot\sum_{k=1}^{n}k^{6}+12\cdot\sum_{k=1}^{n}k^{7}+3\cdot\sum_{k=1}^{n}k^{8}}{24}.
\]

\bigskip{}
At this point we will proceed differently. In the previous calculations
we never adjusted the indices of the sums. Therefore we will start
instead with 
\[
\sum k^{8}=(n+1)\cdot\sum k^{7}-\frac{2\cdot\sum k^{2}-7\cdot\sum k^{4}+14\cdot\sum k^{6}+12\cdot\sum k^{7}+3\cdot\sum k^{8}}{24}.
\]
 We rewrite it as 
\[
\sum k^{8}=(n+1)\cdot\sum k^{7}-\frac{1}{12}\cdot\sum k^{2}+\frac{7}{24}\cdot\sum k^{4}-\frac{7}{12}\cdot\sum k^{6}-\frac{1}{2}\cdot\sum k^{7}-\frac{1}{8}\cdot\sum k^{8}
\]
\[
\frac{9}{8}\cdot\sum k^{8}=\frac{2n+1}{2}\cdot\sum k^{7}-\frac{1}{12}\cdot\sum k^{2}+\frac{7}{24}\cdot\sum k^{4}-\frac{7}{12}\cdot\sum k^{6}.
\]

In order to simplify the calculations we want to delay introducing
the coefficients for the terms. In Section \ref{sec:Emerging-Patterns}
we introduced the notation $E_{2m}(n)$ and $O_{2m+1}(n)$. Let us
simplify it to $E_{2m}$ and $O_{2m+1}$. Then we may rewrite the
last expression as 
\begin{eqnarray*}
\frac{9}{8}\cdot\sum k^{8} & = & \frac{2n+1}{2}\cdot O_{7}\cdot\sum k^{3}\\
 &  & -\frac{1}{12}\cdot\sum k^{2}+\frac{7}{24}\cdot E_{4}\cdot\sum k^{2}-\frac{7}{12}\cdot E_{6}\cdot\sum k^{2},
\end{eqnarray*}
 where $O_{7},E_{4},$ and $E_{6}$ stem from the previous expressions
for $\sum k^{7},\sum k^{4},$ and $\sum k^{6}$.

Since we believe the final expression will be of the form $\sum k^{8}=E_{8}\cdot\sum k^{2}$,
we rewrite the term on the left side as follows: 
\[
\frac{2n+1}{2}\cdot O_{7}\cdot\frac{n\left(n+1\right)}{2}\cdot\frac{n\left(n+1\right)}{2}=\frac{3}{3}\cdot\frac{2n+1}{2}\cdot O_{7}\cdot\frac{n\left(n+1\right)}{2}\cdot\frac{n\left(n+1\right)}{2}
\]
\[
=\frac{3}{2}\cdot\frac{n\left(n+1\right)}{2}\cdot O_{7}\cdot\frac{2n+1}{3}\cdot\frac{n\left(n+1\right)}{2}
\]
\[
=\frac{3}{2}\cdot\frac{n(n+1)}{2}\cdot O_{7}\cdot\sum k^{2}.
\]
 In total we have 
\begin{eqnarray}
\frac{9}{8}\cdot\sum k^{8} & = & \frac{3}{2}\cdot\frac{n(n+1)}{2}\cdot O_{7}\cdot\sum k^{2}\nonumber \\
 &  & -\frac{1}{12}\cdot\sum k^{2}+\frac{7}{24}\cdot E_{4}\cdot\sum k^{2}-\frac{7}{12}\cdot E_{6}\cdot\sum k^{2}.\label{eq:17}
\end{eqnarray}
 Now we are ready to substitute the true coefficients. Only, we are
going to use a different set of coefficients.

\bigskip{}
If we look back to the summary of Part 1 in Section \ref{sec:Summary},
we will remind ourselves we have been carrying two forms of the expressions.
Now it will be advantageous to return to the soft forms: 
\[
\sum k^{4}=\frac{1}{5}\cdot\left(6\cdot\frac{n\left(n+1\right)}{2}-1\right)\cdot\sum k^{2}
\]
\[
\sum k^{6}=\frac{1}{7}\cdot\left(12\cdot\left(\frac{n\left(n+1\right)}{2}\right)^{2}-6\cdot\frac{n\left(n+1\right)}{2}+1\right)\cdot\sum k^{2}
\]
\[
\sum k^{7}=\frac{1}{3}\cdot\left(6\cdot\left(\frac{n(n+1)}{2}\right)^{2}-4\cdot\frac{n(n+1)}{2}+1\right)\cdot\sum k^{3}.
\]
 We substitute the coefficients into expression \ref{eq:17} to get
\begin{eqnarray*}
\frac{9}{8}\cdot\sum k^{8} & = & \frac{3}{2}\cdot\frac{n(n+1)}{2}\cdot\frac{1}{3}\cdot\left(6\cdot\left(\frac{n(n+1)}{2}\right)^{2}-4\cdot\frac{n(n+1)}{2}+1\right)\cdot\sum k^{2}\\
 &  & -\frac{1}{12}\cdot\sum k^{2}\\
 &  & +\frac{7}{24}\cdot\frac{1}{5}\cdot\left(6\cdot\frac{n\left(n+1\right)}{2}-1\right)\cdot\sum k^{2}\\
 &  & -\frac{7}{12}\cdot\frac{1}{7}\cdot\left(12\cdot\left(\frac{n\left(n+1\right)}{2}\right)^{2}-6\cdot\frac{n\left(n+1\right)}{2}+1\right)\cdot\sum k^{2}
\end{eqnarray*}
\begin{eqnarray*}
\frac{9}{8}\cdot\sum k^{8} & = & \left(3\cdot\left(\frac{n(n+1)}{2}\right)^{3}-2\cdot\left(\frac{n(n+1)}{2}\right)^{2}+\frac{1}{2}\cdot\frac{n(n+1)}{2}\right)\cdot\sum k^{2}\\
 &  & -\frac{1}{12}\cdot\sum k^{2}\\
 &  & +\left(\frac{42}{120}\cdot\frac{n\left(n+1\right)}{2}-\frac{7}{120}\right)\cdot\sum k^{2}\\
 &  & +\left(-\left(\frac{n\left(n+1\right)}{2}\right)^{2}+\frac{1}{2}\cdot\frac{n\left(n+1\right)}{2}-\frac{1}{12}\right)\cdot\sum k^{2}
\end{eqnarray*}
\begin{eqnarray*}
\frac{9}{8}\cdot\sum k^{8} & = & 3\cdot\left(\frac{n(n+1)}{2}\right)^{3}\cdot\sum k^{2}\\
 &  & -3\cdot\left(\frac{n(n+1)}{2}\right)^{2}\cdot\sum k^{2}\\
 &  & +\left(\frac{1}{2}+\frac{42}{120}+\frac{1}{2}\right)\cdot\frac{n\left(n+1\right)}{2}\cdot\sum k^{2}\\
 &  & +\left(-\frac{1}{12}-\frac{7}{120}-\frac{1}{12}\right)\cdot\sum k^{2}
\end{eqnarray*}
\begin{eqnarray*}
\frac{9}{8}\cdot\sum k^{8} & = & 3\cdot\left(\frac{n(n+1)}{2}\right)^{3}\cdot\sum k^{2}\\
 &  & -3\cdot\left(\frac{n(n+1)}{2}\right)^{2}\cdot\sum k^{2}\\
 &  & +\frac{27}{20}\cdot\frac{n(n+1)}{2}\cdot\sum k^{2}\\
 &  & -\frac{9}{40}\cdot\sum k^{2},
\end{eqnarray*}
 which we rewrite as 
\begin{eqnarray*}
\sum k^{8} & = & \frac{8}{9}\cdot\left(3\cdot\left(\frac{n(n+1)}{2}\right)^{3}-3\cdot\left(\frac{n(n+1)}{2}\right)^{2}+\frac{27}{20}\cdot\frac{n(n+1)}{2}-\frac{9}{40}\right)\cdot\sum k^{2}\\
 &  & =\left(\frac{8}{3}\cdot\left(\frac{n(n+1)}{2}\right)^{3}-\frac{8}{3}\cdot\left(\frac{n(n+1)}{2}\right)^{2}+\frac{6}{5}\cdot\frac{n(n+1)}{2}-\frac{1}{5}\right)\cdot\sum k^{2}.
\end{eqnarray*}
 Another way to write it is 
\begin{equation}
\sum k^{8}=\frac{1}{9}\cdot\left(24\cdot\left(\frac{n(n+1)}{2}\right)^{3}-24\cdot\left(\frac{n(n+1)}{2}\right)^{2}+\frac{54}{5}\cdot\frac{n(n+1)}{2}-\frac{9}{5}\right)\cdot\sum k^{2}.\label{eq:18}
\end{equation}

We have verified the pattern again.\footnote{\label{fn:Another-way-to}Notice that 
\begin{eqnarray*}
\sum k^{8} & = & \frac{1}{9}\cdot\left(24\cdot\left(\frac{n(n+1)}{2}\right)^{3}-24\cdot\left(\frac{n(n+1)}{2}\right)^{2}+\frac{54}{5}\cdot\frac{n(n+1)}{2}-\frac{9}{5}\right)\cdot\sum k^{2}\\
 &  & =\frac{1}{9}\cdot\left(24\cdot\left(\frac{n(n+1)}{2}\right)^{3}-24\cdot\left(\frac{n(n+1)}{2}\right)^{2}+9\cdot\frac{1}{5}\cdot\left(6\cdot\frac{n(n+1)}{2}-1\right)\right)\cdot\sum k^{2},
\end{eqnarray*}
 which is a curious 
\[
\sum k^{8}=\frac{1}{9}\cdot\left(24\cdot\left(\frac{n(n+1)}{2}\right)^{3}-24\cdot\left(\frac{n(n+1)}{2}\right)^{2}\right)\cdot\sum k^{2}+\sum k^{4}.
\]
} The other form is 
\[
\sum k^{8}=\frac{5(n(n+1))^{3}-10(n(n+1))^{2}+9n(n+1)-3}{15}\cdot\sum k^{2}.
\]
 Both expressions simplify to 
\[
\sum k^{8}=\frac{-3n+20n^{3}-42n^{5}+60n^{7}+45n^{8}+10n^{9}}{90}.
\]
 Are we ready to put all of this together for one, final case?

\bigskip{}

\section{$\sum k^{9}$}

What is 
\[
\sum k^{9}=1^{9}+2^{9}+3^{9}+\cdots+n^{9}\,?
\]
We know that 
\begin{equation}
\sum k^{9}=(n+1)\cdot\sum k^{8}-\sum\sum l^{8}.\label{eq:19}
\end{equation}
 We just derived 
\[
\sum k^{8}=\frac{-3n+20n^{3}-42n^{5}+60n^{7}+45n^{8}+10n^{9}}{90},
\]
 which implies 
\[
\sum\sum l^{8}=\frac{-3\cdot\sum k+20\cdot\sum k^{3}-42\cdot\sum k^{5}+60\cdot\sum k^{7}+45\cdot\sum k^{8}+10\cdot\sum k^{9}}{90}.
\]
 Therefore we may rewrite expression \ref{eq:19} as 
\[
\sum k^{9}=(n+1)\cdot\sum k^{8}-\frac{-3\cdot\sum k+20\cdot\sum k^{3}-42\cdot\sum k^{5}+60\cdot\sum k^{7}+45\cdot\sum k^{8}+10\cdot\sum k^{9}}{90}
\]
\[
\sum k^{9}=(n+1)\cdot\sum k^{8}+\frac{1}{30}\cdot\sum k-\frac{2}{9}\cdot\sum k^{3}+\frac{7}{15}\cdot\sum k^{5}-\frac{2}{3}\cdot\sum k^{7}-\frac{1}{2}\cdot\sum k^{8}-\frac{1}{9}\cdot\sum k^{9}
\]
\[
\frac{10}{9}\cdot\sum k^{9}=\frac{2n+1}{2}\cdot\sum k^{8}+\frac{1}{30}\cdot\sum k-\frac{2}{9}\cdot\sum k^{3}+\frac{7}{15}\cdot\sum k^{5}-\frac{2}{3}\cdot\sum k^{7}
\]
\begin{eqnarray}
\frac{10}{9}\cdot\sum k^{9} & = & \frac{2n+1}{2}\cdot E_{8}\cdot\sum k^{2}\nonumber \\
 &  & +\frac{1}{30}\cdot\sum k\nonumber \\
 &  & -\frac{2}{9}\cdot\sum k^{3}+\frac{7}{15}\cdot O_{5}\cdot\sum k^{3}-\frac{2}{3}\cdot O_{7}\cdot\sum k^{3},\label{eq:20}
\end{eqnarray}
 where 
\[
E_{8}=\frac{8}{3}\cdot\left(\frac{n(n+1)}{2}\right)^{3}-\frac{8}{3}\cdot\left(\frac{n(n+1)}{2}\right)^{2}+\frac{6}{5}\cdot\frac{n(n+1)}{2}-\frac{1}{5},
\]
\[
O_{5}=\frac{1}{3}\cdot\left(4\cdot\frac{n\left(n+1\right)}{2}-1\right),
\]
 and 
\[
O_{7}=\frac{1}{3}\cdot\left(6\cdot\left(\frac{n(n+1)}{2}\right)^{2}-4\cdot\frac{n(n+1)}{2}+1\right).
\]

\bigskip{}
We believe the final expression will be of the form $\sum k^{9}=O_{9}\cdot\sum k^{3}$.
However, analogous to the calculation for $\sum k^{7}$, we will aim
first for $\sum k^{9}=O_{9}\cdot\frac{n(n+1)}{2}\cdot\sum k$. For
expression \ref{eq:20} we rewrite the final terms as 
\[
-\frac{2}{9}\cdot\sum k^{3}+\frac{7}{15}\cdot O_{5}\cdot\sum k^{3}-\frac{2}{3}\cdot O_{7}\cdot\sum k^{3}
\]
\[
=-\frac{2}{9}\cdot\frac{n(n+1)}{2}\cdot\sum k+\frac{7}{15}\cdot O_{5}\cdot\frac{n(n+1)}{2}\cdot\sum k-\frac{2}{3}\cdot\frac{n(n+1)}{2}\cdot O_{7}\cdot\sum k
\]
\begin{eqnarray*}
 & = & -\frac{2}{9}\cdot\frac{n(n+1)}{2}\cdot\sum k+\frac{7}{15}\cdot\frac{1}{3}\cdot\left(4\cdot\frac{n\left(n+1\right)}{2}-1\right)\cdot\frac{n(n+1)}{2}\cdot\sum k\\
 &  & -\frac{2}{3}\cdot\frac{1}{3}\cdot\left(6\cdot\left(\frac{n(n+1)}{2}\right)^{2}-4\cdot\frac{n(n+1)}{2}+1\right)\cdot\frac{n(n+1)}{2}\cdot\sum k
\end{eqnarray*}
\begin{eqnarray*}
 & = & -\frac{2}{9}\cdot\frac{n(n+1)}{2}\cdot\sum k\\
 &  & +\left(\frac{28}{45}\cdot\frac{n(n+1)}{2}-\frac{7}{45}\right)\cdot\frac{n(n+1)}{2}\cdot\sum k\\
 &  & +\left(-\frac{4}{3}\cdot\left(\frac{n(n+1)}{2}\right)^{2}+\frac{8}{9}\cdot\frac{n(n+1)}{2}-\frac{2}{9}\right)\cdot\frac{n(n+1)}{2}\cdot\sum k
\end{eqnarray*}
\begin{eqnarray*}
 & = & -\frac{2}{9}\cdot\frac{n(n+1)}{2}\cdot\sum k\\
 &  & +\left(\frac{28}{45}\cdot\left(\frac{n(n+1)}{2}\right)^{2}-\frac{7}{45}\cdot\frac{n(n+1)}{2}\right)\cdot\sum k\\
 &  & +\left(-\frac{4}{3}\cdot\left(\frac{n(n+1)}{2}\right)^{3}+\frac{8}{9}\cdot\left(\frac{n(n+1)}{2}\right)^{2}-\frac{2}{9}\cdot\frac{n(n+1)}{2}\right)\cdot\sum k
\end{eqnarray*}
\[
=\left(-\frac{4}{3}\cdot\left(\frac{n(n+1)}{2}\right)^{3}+\left(\frac{28}{45}+\frac{8}{9}\right)\cdot\left(\frac{n(n+1)}{2}\right)^{2}+\left(-\frac{2}{9}-\frac{7}{45}-\frac{2}{9}\right)\cdot\left(\frac{n(n+1)}{2}\right)\right)\cdot\sum k.
\]
\begin{equation}
=\left(-\frac{4}{3}\cdot\left(\frac{n(n+1)}{2}\right)^{3}+\frac{68}{45}\cdot\left(\frac{n(n+1)}{2}\right)^{2}-\frac{3}{5}\cdot\left(\frac{n(n+1)}{2}\right)\right)\cdot\sum k.\label{eq:21}
\end{equation}

\bigskip{}
Next, we rewrite the first term as 
\begin{eqnarray*}
\frac{2n+1}{2}\cdot E_{8}\cdot\sum k^{2} & = & \frac{2n+1}{2}\cdot E_{8}\cdot\frac{2n+1}{3}\cdot\sum k\\
 &  & =\frac{(2n+1)^{2}}{6}\cdot E_{8}\cdot\sum k\\
 &  & =\left(\frac{4}{3}\cdot\frac{n(n+1)}{2}+\frac{1}{6}\right)\cdot E_{8}\cdot\sum k
\end{eqnarray*}
 and then substitute the expression for $E_{8}$: 
\[
\left(\frac{4}{3}\cdot\frac{n(n+1)}{2}+\frac{1}{6}\right)\cdot\left(\frac{8}{3}\cdot\left(\frac{n(n+1)}{2}\right)^{3}-\frac{8}{3}\cdot\left(\frac{n(n+1)}{2}\right)^{2}+\frac{6}{5}\cdot\frac{n(n+1)}{2}-\frac{1}{5}\right)\cdot\sum k
\]
\begin{eqnarray*}
 & = & \left(\frac{32}{9}\cdot\left(\frac{n(n+1)}{2}\right)^{4}-\frac{32}{9}\cdot\left(\frac{n(n+1)}{2}\right)^{3}+\frac{24}{15}\cdot\left(\frac{n(n+1)}{2}\right)^{2}-\frac{4}{15}\cdot\frac{n(n+1)}{2}\right)\cdot\sum k\\
 &  & +\left(\frac{8}{18}\cdot\left(\frac{n(n+1)}{2}\right)^{3}-\frac{8}{18}\cdot\left(\frac{n(n+1)}{2}\right)^{2}+\frac{1}{5}\cdot\frac{n(n+1)}{2}-\frac{1}{30}\right)\cdot\sum k
\end{eqnarray*}
\begin{eqnarray*}
 & = & \frac{32}{9}\cdot\left(\frac{n(n+1)}{2}\right)^{4}\cdot\sum k\\
 &  & +\left(-\frac{32}{9}+\frac{8}{18}\right)\cdot\left(\frac{n(n+1)}{2}\right)^{3}\cdot\sum k\\
 &  & +\left(\frac{24}{15}-\frac{8}{18}\right)\cdot\left(\frac{n(n+1)}{2}\right)^{2}\cdot\sum k\\
 &  & +\left(-\frac{4}{15}+\frac{1}{5}\right)\cdot\frac{n(n+1)}{2}\cdot\sum k\\
 &  & -\frac{1}{30}\cdot\sum k
\end{eqnarray*}
\begin{equation}
=\left(\frac{32}{9}\cdot\left(\frac{n(n+1)}{2}\right)^{4}-\frac{28}{9}\cdot\left(\frac{n(n+1)}{2}\right)^{3}+\frac{52}{45}\cdot\left(\frac{n(n+1)}{2}\right)^{2}-\frac{1}{15}\cdot\frac{n(n+1)}{2}\right)\cdot\sum k-\frac{1}{30}\cdot\sum k.\label{eq:22}
\end{equation}
 The lone term of $-\frac{1}{30}\cdot\sum k$ should make us happy.
It cancels with the $\frac{1}{30}\cdot\sum k$ in expression \ref{eq:20}.

\bigskip{}
Finally, if we put together expressions \ref{eq:20}, \ref{eq:21},
and \ref{eq:22} then we get 
\begin{eqnarray*}
\frac{10}{9}\cdot\sum k^{9} & = & \left(\frac{32}{9}\cdot\left(\frac{n(n+1)}{2}\right)^{4}-\frac{28}{9}\cdot\left(\frac{n(n+1)}{2}\right)^{3}\right)\cdot\sum k\\
 &  & +\left(\frac{52}{45}\cdot\left(\frac{n(n+1)}{2}\right)^{2}-\frac{1}{15}\cdot\frac{n(n+1)}{2}\right)\cdot\sum k\\
 &  & -\frac{1}{30}\cdot\sum k+\frac{1}{30}\cdot\sum k\\
 &  & +\left(-\frac{4}{3}\cdot\left(\frac{n(n+1)}{2}\right)^{3}+\frac{68}{45}\cdot\left(\frac{n(n+1)}{2}\right)^{2}-\frac{3}{5}\cdot\left(\frac{n(n+1)}{2}\right)\right)\cdot\sum k
\end{eqnarray*}
\[
\frac{10}{9}\cdot\sum k^{9}=\left(\frac{32}{9}\cdot\left(\frac{n(n+1)}{2}\right)^{4}-\frac{40}{9}\cdot\left(\frac{n(n+1)}{2}\right)^{3}+\frac{8}{3}\cdot\left(\frac{n(n+1)}{2}\right)^{2}-\frac{2}{3}\cdot\frac{n(n+1)}{2}\right)\cdot\sum k
\]
\begin{eqnarray}
\sum k^{9} & = & \left(\frac{16}{5}\cdot\left(\frac{n(n+1)}{2}\right)^{4}-4\cdot\left(\frac{n(n+1)}{2}\right)^{3}+\frac{12}{5}\cdot\left(\frac{n(n+1)}{2}\right)^{2}-\frac{3}{5}\cdot\frac{n(n+1)}{2}\right)\cdot\sum k\nonumber \\
 &  & =\left(\frac{16}{5}\cdot\left(\frac{n(n+1)}{2}\right)^{3}-4\cdot\left(\frac{n(n+1)}{2}\right)^{2}+\frac{12}{5}\cdot\frac{n(n+1)}{2}-\frac{3}{5}\right)\cdot\sum k^{3}\nonumber \\
 &  & =\frac{1}{5}\cdot\left(16\cdot\left(\frac{n(n+1)}{2}\right)^{3}-20\cdot\left(\frac{n(n+1)}{2}\right)^{2}+12\cdot\frac{n(n+1)}{2}-3\right)\cdot\sum k^{3}.\label{eq:23}
\end{eqnarray}
 The other form is 
\[
\sum k^{9}=\frac{2(n(n+1))^{3}-5(n(n+1))^{2}+6n(n+1)-3}{5}\cdot\sum k^{3}.
\]
 Both expressions simplify to 
\[
\sum k^{9}=\frac{-3n^{2}+10n^{4}-14n^{6}+15n^{8}+10n^{9}+2n^{10}}{20}.
\]
 We got our $\sum k^{9}=O_{9}\cdot\sum k^{3}$.\pagebreak{}

\section{Summary of Part 2}

Building upon the methods and results of Part 1, we guessed the sums
might follow patterns for even and odd powers: 
\[
\sum_{k=1}^{n}k^{2m}=E_{2m}(n)\cdot\sum_{k=1}^{n}k^{2}
\]
\[
\sum_{k=1}^{n}k^{2m+1}=O_{2m+1}(n)\cdot\sum_{k=1}^{n}k^{3},
\]
 where $E_{2m}(n)$ and $O_{2m+1}(n)$ were rational expressions involving
$n$. We were able to verify them for the next cases of $m=6,7$:
\begin{eqnarray*}
\sum k^{6} & = & \frac{1}{7}\cdot\left(12\cdot\left(\frac{n\left(n+1\right)}{2}\right)^{2}-6\cdot\frac{n\left(n+1\right)}{2}+1\right)\cdot\sum k^{2}\\
 &  & =\left(\frac{3\left(n\left(n+1\right)\right)^{2}-3n\left(n+1\right)+1}{7}\right)\cdot\sum k^{2}
\end{eqnarray*}
\begin{eqnarray*}
\sum k^{7} & = & \frac{1}{3}\cdot\left(6\cdot\left(\frac{n(n+1)}{2}\right)^{2}-4\cdot\frac{n(n+1)}{2}+1\right)\cdot\sum k^{3}\\
 &  & =\left(\frac{3(n(n+1))^{2}-4n(n+1)+2}{6}\right)\cdot\sum k^{3}.
\end{eqnarray*}
We suspected the same would be true for $m=8,9$, but the calculations
became too difficult. We modified our approach and derived 
\begin{eqnarray*}
\sum k^{8} & = & \frac{1}{9}\cdot\left(24\cdot\left(\frac{n(n+1)}{2}\right)^{3}-24\cdot\left(\frac{n(n+1)}{2}\right)^{2}+\frac{54}{5}\cdot\frac{n(n+1)}{2}-\frac{9}{5}\right)\cdot\sum k^{2}\\
 &  & =\frac{5(n(n+1))^{3}-10(n(n+1))^{2}+9n(n+1)-3}{15}\cdot\sum k^{2}
\end{eqnarray*}
\begin{eqnarray*}
\sum k^{9} & = & \frac{1}{5}\cdot\left(16\cdot\left(\frac{n(n+1)}{2}\right)^{3}-20\cdot\left(\frac{n(n+1)}{2}\right)^{2}+12\cdot\frac{n(n+1)}{2}-3\right)\cdot\sum k^{3}\\
 &  & =\frac{2(n(n+1))^{3}-5(n(n+1))^{2}+6n(n+1)-3}{5}\cdot\sum k^{3},
\end{eqnarray*}
 which verified the patterns again. Are we able to guess general expressions
for $E_{2m}$ and $O_{2m+1}$?

\pagebreak{}

\part*{Part 3}

\section{A Long List\label{sec:A-Long-List}}

We have assembled a long list of results. Let us look at all of them
together: 
\begin{eqnarray*}
\\
\sum k^{2} & = & \frac{2n+1}{3}\cdot\frac{n\left(n+1\right)}{2}\\
\sum k^{3} & = & \left(\frac{n\left(n+1\right)}{2}\right)^{2}\\
\sum k^{4} & = & \frac{1}{5}\cdot\left(6\cdot\frac{n\left(n+1\right)}{2}-1\right)\cdot\sum k^{2}\\
\sum k^{5} & = & \frac{1}{3}\cdot\left(4\cdot\frac{n\left(n+1\right)}{2}-1\right)\cdot\sum k^{3}\\
\sum k^{6} & = & \frac{1}{7}\cdot\left(12\cdot\left(\frac{n\left(n+1\right)}{2}\right)^{2}-6\cdot\frac{n\left(n+1\right)}{2}+1\right)\cdot\sum k^{2}\\
\sum k^{7} & = & \frac{1}{3}\cdot\left(6\cdot\left(\frac{n(n+1)}{2}\right)^{2}-4\cdot\frac{n(n+1)}{2}+1\right)\cdot\sum k^{3}\\
\sum k^{8} & = & \frac{1}{9}\cdot\left(24\cdot\left(\frac{n(n+1)}{2}\right)^{3}-24\cdot\left(\frac{n(n+1)}{2}\right)^{2}+\frac{54}{5}\cdot\frac{n(n+1)}{2}-\frac{9}{5}\right)\cdot\sum k^{2}\\
\sum k^{9} & = & \frac{1}{5}\cdot\left(16\cdot\left(\frac{n(n+1)}{2}\right)^{3}-20\cdot\left(\frac{n(n+1)}{2}\right)^{2}+12\cdot\frac{n(n+1)}{2}-3\right)\cdot\sum k^{3},\\
\end{eqnarray*}
 and of course, $\sum k=\frac{n(n+1)}{2}.$ Let us group them in even
and odd powers: 
\begin{eqnarray*}
\\
\sum k^{2} & = & \frac{2n+1}{3}\cdot\frac{n\left(n+1\right)}{2}\\
\sum k^{4} & = & \frac{1}{5}\cdot\left(6\cdot\frac{n\left(n+1\right)}{2}-1\right)\cdot\sum k^{2}\\
\sum k^{6} & = & \frac{1}{7}\cdot\left(12\cdot\left(\frac{n\left(n+1\right)}{2}\right)^{2}-6\cdot\frac{n\left(n+1\right)}{2}+1\right)\cdot\sum k^{2}\\
\sum k^{8} & = & \frac{1}{9}\cdot\left(24\cdot\left(\frac{n(n+1)}{2}\right)^{3}-24\cdot\left(\frac{n(n+1)}{2}\right)^{2}+\frac{54}{5}\cdot\frac{n(n+1)}{2}-\frac{9}{5}\right)\cdot\sum k^{2}
\end{eqnarray*}
\pagebreak{}
\begin{eqnarray*}
\sum k^{3} & = & \left(\frac{n\left(n+1\right)}{2}\right)^{2}\\
\sum k^{5} & = & \frac{1}{3}\cdot\left(4\cdot\frac{n\left(n+1\right)}{2}-1\right)\cdot\sum k^{3}\\
\sum k^{7} & = & \frac{1}{3}\cdot\left(6\cdot\left(\frac{n(n+1)}{2}\right)^{2}-4\cdot\frac{n(n+1)}{2}+1\right)\cdot\sum k^{3}\\
\sum k^{9} & = & \frac{1}{5}\cdot\left(16\cdot\left(\frac{n(n+1)}{2}\right)^{3}-20\cdot\left(\frac{n(n+1)}{2}\right)^{2}+12\cdot\frac{n(n+1)}{2}-3\right)\cdot\sum k^{3}.\\
\end{eqnarray*}
For the odd powers, something is amiss with the leading fractions.
Let us try 
\begin{eqnarray*}
\\
\sum k^{2} & = & \frac{2n+1}{3}\cdot\frac{n\left(n+1\right)}{2}\\
\sum k^{4} & = & \frac{1}{5}\cdot\left(6\cdot\frac{n\left(n+1\right)}{2}-1\right)\cdot\sum k^{2}\\
\sum k^{6} & = & \frac{1}{7}\cdot\left(12\cdot\left(\frac{n\left(n+1\right)}{2}\right)^{2}-6\cdot\frac{n\left(n+1\right)}{2}+1\right)\cdot\sum k^{2}\\
\sum k^{8} & = & \frac{1}{9}\cdot\left(24\cdot\left(\frac{n(n+1)}{2}\right)^{3}-24\cdot\left(\frac{n(n+1)}{2}\right)^{2}+\frac{54}{5}\cdot\frac{n(n+1)}{2}-\frac{9}{5}\right)\cdot\sum k^{2}
\end{eqnarray*}
\begin{eqnarray*}
\sum k^{3} & = & \left(\frac{n\left(n+1\right)}{2}\right)^{2}\\
\sum k^{5} & = & \frac{1}{6}\cdot\left(8\cdot\frac{n\left(n+1\right)}{2}-2\right)\cdot\sum k^{3}\\
\sum k^{7} & = & \frac{1}{8}\cdot\left(\frac{1}{3}\cdot\left(48\cdot\left(\frac{n(n+1)}{2}\right)^{2}-32\cdot\frac{n(n+1)}{2}+8\right)\right)\cdot\sum k^{3}\\
\sum k^{9} & = & \frac{1}{10}\cdot\left(32\cdot\left(\frac{n(n+1)}{2}\right)^{3}-40\cdot\left(\frac{n(n+1)}{2}\right)^{2}+24\cdot\frac{n(n+1)}{2}-6\right)\cdot\sum k^{3}.\\
\end{eqnarray*}
 Yes, that looks better. Do we notice any patterns?

\bigskip{}

We observe the following:
\begin{enumerate}
\item the leading fractions for the sums seem to be based on the powers
of the terms: for example, $\frac{1}{5}$ for $\sum k^{4}$ and $\frac{1}{10}$
for $\sum k^{9}$. We expect to find $\frac{1}{11}$ for $\sum k^{10}$
and $\frac{1}{12}$ for $\sum k^{11}$.
\item the terms of $\frac{n(n+1)}{2}$ in the coefficients for the sums
are raised to the powers 
\begin{eqnarray*}
 & 1\\
 & 2,1\\
 & 3,2,1.
\end{eqnarray*}
 We believe the next cases will be 
\begin{eqnarray*}
 & 4,3,2,1,\\
 & 5,4,3,2,1,
\end{eqnarray*}
 and so forth.
\item as for the coefficients for the terms of $\frac{n(n+1)}{2}$ themselves,
either they do not follow a pattern or they follow one which remains
a mystery still. We have 
\begin{eqnarray*}
 & 6,-1\\
 & 12,-6,1\\
 & 24,-24,\frac{54}{5},-\frac{9}{5}
\end{eqnarray*}
 and 
\begin{eqnarray*}
 & 8,-2\\
 & 16,-\frac{32}{3},\frac{8}{3}\\
 & 32,-40,24,-6.
\end{eqnarray*}
 Why is a multiple of $\frac{1}{5}$ introduced? of $\frac{1}{3}$?
We do not know. If we clear the fractions then we get 
\begin{eqnarray*}
 & 6,-1\\
 & 12,-6,1\\
 & 120,-120,54,-9
\end{eqnarray*}
 and 
\begin{eqnarray*}
 & 8,-2\\
 & 48,-32,8\\
 & 32,-40,24,-6.
\end{eqnarray*}
 It still is hard to discern a pattern. About the only thing we can
pinpoint are the alternating signs in front of the coefficients.
\end{enumerate}
Honestly, we are disappointed. After all of this work we expected
the general patterns for $\sum k^{2m}$ and $\sum k^{2m+1}$ to fall
into our laps. We are a far cry from that. Worse, we remember how
lengthy the calculations were for $\sum k^{8}$ and $\sum k^{9}$
and are reluctant to attempt them for $\sum k^{10}$ and $\sum k^{11}$.
What else can we do?

\section{A Different Point of View\label{sec:A-Different-Point-of-View}}

Perhaps we are looking at the expressions too directly. Perhaps we
would benefit from taking a different point of view. Let us summarize
our results as follows:
\begin{enumerate}
\item we have observed a pattern for even powers and a pattern for odd powers.
The patterns have a close resemblance to one another.
\item after distinguishing the two patterns, the forms of the total sums
have become less important and the forms of the coefficients for the
sums have become more important.
\end{enumerate}
Let us turn our attention to the coefficients.

\bigskip{}
For the calculation for $\sum k^{8}$ in Section \ref{sec:Refining-Our-Approach},
in Footnote \ref{fn:Another-way-to} we remarked we noticed something
curious: 
\begin{eqnarray*}
\sum k^{8} & = & \frac{1}{9}\cdot\left(24\cdot\left(\frac{n(n+1)}{2}\right)^{3}-24\cdot\left(\frac{n(n+1)}{2}\right)^{2}+\frac{54}{5}\cdot\frac{n(n+1)}{2}-\frac{9}{5}\right)\cdot\sum k^{2}\\
 &  & =\frac{1}{9}\cdot\left(24\cdot\left(\frac{n(n+1)}{2}\right)^{3}-24\cdot\left(\frac{n(n+1)}{2}\right)^{2}+9\cdot\frac{1}{5}\cdot\left(6\cdot\frac{n(n+1)}{2}-1\right)\right)\cdot\sum k^{2}\\
 &  & =\frac{1}{9}\cdot\left(24\cdot\left(\frac{n(n+1)}{2}\right)^{3}-24\cdot\left(\frac{n(n+1)}{2}\right)^{2}\right)\cdot\sum k^{2}+\sum k^{4},
\end{eqnarray*}
 where 
\[
\sum k^{4}=\frac{1}{5}\cdot\left(6\cdot\frac{n\left(n+1\right)}{2}-1\right)\cdot\sum k^{2}.
\]
 The expression for $\sum k^{4}$ ended up in the expression for $\sum k^{8}$.
How? More curious, if we look at 
\[
\sum k^{6}=\frac{1}{7}\cdot\left(12\cdot\left(\frac{n\left(n+1\right)}{2}\right)^{2}-6\cdot\frac{n\left(n+1\right)}{2}+1\right)\cdot\sum k^{2}
\]
 then we notice 
\begin{eqnarray*}
\sum k^{6} & = & \frac{1}{7}\cdot\left(12\cdot\left(\frac{n\left(n+1\right)}{2}\right)^{2}-5\cdot\frac{1}{5}\left(6\cdot\frac{n\left(n+1\right)}{2}-1\right)\right)\cdot\sum k^{2}\\
 &  & =\frac{12}{7}\cdot\left(\frac{n\left(n+1\right)}{2}\right)^{2}\cdot\sum k^{2}-\frac{5}{7}\cdot\sum k^{4}.
\end{eqnarray*}
 $\sum k^{4}$ appeared again. What is the explanation?

\subsection{Even}

In Section \ref{sec:Emerging-Patterns} we introduced the notation
$E_{2m}(n)$ and $O_{2m+1}(n)$, which refer to the coefficients for
$\sum k^{2m}$ and $\sum k^{2m+1}$, respectively. Let us return to
it.

\bigskip{}
For $\sum k^{4},\sum k^{6},$ and $\sum k^{8}$ we look at the coefficients
rather than the full expressions: 
\[
E_{4}=\frac{1}{5}\cdot\left(6\cdot\frac{n\left(n+1\right)}{2}-1\right)
\]
\[
E_{6}=\frac{1}{7}\cdot\left(12\cdot\left(\frac{n\left(n+1\right)}{2}\right)^{2}-6\cdot\frac{n\left(n+1\right)}{2}+1\right)
\]
\[
E_{8}=\frac{1}{9}\cdot\left(24\cdot\left(\frac{n(n+1)}{2}\right)^{3}-24\cdot\left(\frac{n(n+1)}{2}\right)^{2}+\frac{9}{5}\cdot\left(6\cdot\frac{n(n+1)}{2}-1\right)\right).
\]
 The new observations allow us to write 
\[
E_{6}=\frac{12}{7}\cdot\left(\frac{n\left(n+1\right)}{2}\right)^{2}-\frac{5}{7}\cdot E_{4}
\]
\[
E_{8}=\frac{1}{9}\cdot\left(24\cdot\left(\frac{n(n+1)}{2}\right)^{3}-24\cdot\left(\frac{n(n+1)}{2}\right)^{2}\right)+E_{4}.
\]
 What about $E_{4}$? In $E_{6}$ and $E_{8}$ we noticed the appearance
of $E_{4}$. In 
\[
E_{4}=\frac{1}{5}\cdot\left(6\cdot\frac{n\left(n+1\right)}{2}-1\right)
\]
 do we notice anything? What is the term of $-1$?

\bigskip{}
Ah, we have forgotten something. Remember that 
\[
\sum k^{2}=\frac{2n+1}{3}\cdot\frac{n\left(n+1\right)}{2}.
\]
 If we write 
\[
\sum k^{2}=\frac{2n+1}{3}\cdot\frac{n\left(n+1\right)}{2}=1\cdot\frac{2n+1}{3}\cdot\frac{n\left(n+1\right)}{2}=1\cdot\sum k^{2}
\]
 then we get $E_{2}=1$. This allows us to rewrite $E_{4}$ as 
\[
E_{4}=\frac{1}{5}\cdot\left(6\cdot\frac{n\left(n+1\right)}{2}-E_{2}\right)
\]
\[
5E_{4}=6\cdot\frac{n\left(n+1\right)}{2}-E_{2}
\]
 or 
\begin{equation}
E_{2}+5E_{4}=6\cdot\frac{n\left(n+1\right)}{2}.\label{eq:24}
\end{equation}
 In an analogous fashion we may rewrite $E_{6}$ as 
\[
E_{6}=\frac{12}{7}\cdot\left(\frac{n\left(n+1\right)}{2}\right)^{2}-\frac{5}{7}\cdot E_{4}
\]
\[
7E_{6}=12\cdot\left(\frac{n\left(n+1\right)}{2}\right)^{2}-5E_{4}
\]
 or 
\begin{equation}
5E_{4}+7E_{6}=12\cdot\left(\frac{n\left(n+1\right)}{2}\right)^{2},\label{eq:25}
\end{equation}
 and $E_{8}$ as 
\[
E_{8}=\frac{1}{9}\cdot\left(24\cdot\left(\frac{n(n+1)}{2}\right)^{3}-24\cdot\left(\frac{n(n+1)}{2}\right)^{2}\right)+E_{4}
\]
\[
-E_{4}+E_{8}=\frac{1}{9}\cdot\left(24\cdot\left(\frac{n(n+1)}{2}\right)^{3}-24\cdot\left(\frac{n(n+1)}{2}\right)^{2}\right)
\]
\[
-9E_{4}+9E_{8}=24\cdot\left(\frac{n(n+1)}{2}\right)^{3}-2\cdot(5E_{4}+7E_{6})
\]
 or 
\begin{equation}
E_{4}+14E_{6}+9E_{8}=24\cdot\left(\frac{n(n+1)}{2}\right)^{3}.\label{eq:26}
\end{equation}
 Together we have 
\begin{eqnarray}
 & E_{2}=1\nonumber \\
 & E_{2}+5E_{4}=6\cdot\frac{n\left(n+1\right)}{2}\nonumber \\
 & 5E_{4}+7E_{6}=12\cdot\left(\frac{n\left(n+1\right)}{2}\right)^{2}\nonumber \\
 & E_{4}+14E_{6}+9E_{8}=24\cdot\left(\frac{n(n+1)}{2}\right)^{3}.\label{eq:27}
\end{eqnarray}

This is remarkable. On the right side of the expressions we have 
\[
6\cdot\frac{n\left(n+1\right)}{2},\,12\cdot\left(\frac{n\left(n+1\right)}{2}\right)^{2},\,24\cdot\left(\frac{n\left(n+1\right)}{2}\right)^{3}.
\]
 We see the terms of $\frac{n\left(n+1\right)}{2}$ are raised to
the powers $1,2,3$ and we notice the simple pattern of 
\[
6=2\cdot3,\,12=2^{2}\cdot3,\,24=2^{3}\cdot3.
\]
 On the left side of the expressions we have sums involving $E_{2},\,E_{4},\,E_{6},$
and $E_{8}$ and coefficients in integers. We might not be able to
explain the integers yet, but we are confident enough to conjecture
the next expression will be 
\begin{equation}
e_{4}E_{4}+e_{6}E_{6}+e_{8}E_{8}+e_{10}E_{10}=48\cdot\left(\frac{n(n+1)}{2}\right)^{4}\label{eq:28}
\end{equation}
 for some integers $e_{4},e_{6},e_{8},e_{10}$. Do the coefficients
for sums of odd powers follow an analogous pattern?

\subsection{Odd}

We remind ourselves that 
\[
O_{3}=1
\]
\begin{eqnarray*}
O_{5} & = & \frac{1}{6}\cdot\left(8\cdot\frac{n\left(n+1\right)}{2}-2\right)\\
 &  & =\frac{1}{3}\cdot\left(4\cdot\frac{n\left(n+1\right)}{2}-1\right)
\end{eqnarray*}
\begin{eqnarray*}
O_{7} & = & \frac{1}{8}\cdot\left(\frac{1}{3}\cdot\left(48\cdot\left(\frac{n(n+1)}{2}\right)^{2}-32\cdot\frac{n(n+1)}{2}+8\right)\right)\\
 &  & =\frac{1}{3}\cdot\left(6\cdot\left(\frac{n(n+1)}{2}\right)^{2}-4\cdot\frac{n(n+1)}{2}+1\right)
\end{eqnarray*}
\[
O_{9}=\frac{1}{10}\cdot\left(32\cdot\left(\frac{n(n+1)}{2}\right)^{3}-40\cdot\left(\frac{n(n+1)}{2}\right)^{2}+24\cdot\frac{n(n+1)}{2}-6\right).
\]
 We rewrite $O_{5}$ as 
\[
O_{5}=\frac{1}{3}\cdot\left(4\cdot\frac{n\left(n+1\right)}{2}-O_{3}\right)
\]
\[
3O_{5}=4\cdot\frac{n\left(n+1\right)}{2}-O_{3}
\]
 or 
\begin{equation}
O_{3}+3O_{5}=4\cdot\frac{n\left(n+1\right)}{2}.\label{eq:29}
\end{equation}
 We rewrite $O_{7}$ as 
\[
O_{7}=\frac{1}{3}\cdot\left(6\cdot\left(\frac{n(n+1)}{2}\right)^{2}-3\cdot\frac{1}{3}\cdot\left(4\cdot\frac{n(n+1)}{2}-1\right)\right)
\]
\[
3O_{7}=6\cdot\left(\frac{n(n+1)}{2}\right)^{2}-3O_{5}
\]
\[
3O_{5}+3O_{7}=6\cdot\left(\frac{n(n+1)}{2}\right)^{2}
\]
\[
O_{5}+O_{7}=2\cdot\left(\frac{n(n+1)}{2}\right)^{2}
\]
 or 
\begin{equation}
4O_{5}+4O_{7}=8\cdot\left(\frac{n(n+1)}{2}\right)^{2}.\label{eq:30}
\end{equation}
 We rewrite $O_{9}$ as 
\[
O_{9}=\frac{1}{10}\cdot\left(32\cdot\left(\frac{n(n+1)}{2}\right)^{3}-40\cdot\left(\frac{n(n+1)}{2}\right)^{2}+18\cdot\frac{1}{3}\cdot\left(4\cdot\frac{n(n+1)}{2}-1\right)\right)
\]
\[
10O_{9}=32\cdot\left(\frac{n(n+1)}{2}\right)^{3}-5\cdot(4O_{5}+4O_{7})+18\cdot O_{5}
\]
\[
10O_{9}=32\cdot\left(\frac{n(n+1)}{2}\right)^{3}-2O_{5}-20O_{7}
\]
\[
2O_{5}+20O_{7}+10O_{9}=32\cdot\left(\frac{n(n+1)}{2}\right)^{3}
\]
 or 
\begin{equation}
O_{5}+10\cdot O_{7}+5O_{9}=16\cdot\left(\frac{n(n+1)}{2}\right)^{3}.\label{eq:31}
\end{equation}
 Together we have 
\begin{eqnarray}
 & O_{3}=1\nonumber \\
 & O_{3}+3O_{5}=4\cdot\frac{n\left(n+1\right)}{2}\nonumber \\
 & 4O_{5}+4O_{7}=8\cdot\left(\frac{n(n+1)}{2}\right)^{2}\nonumber \\
 & O_{5}+10\cdot O_{7}+5O_{9}=16\cdot\left(\frac{n(n+1)}{2}\right)^{3}.\label{eq:32}
\end{eqnarray}

Again, this is remarkable. On the right side we have 
\[
4\cdot\frac{n\left(n+1\right)}{2},\,8\cdot\left(\frac{n\left(n+1\right)}{2}\right)^{2},\,16\cdot\left(\frac{n\left(n+1\right)}{2}\right)^{3}.
\]
 We see the terms of $\frac{n\left(n+1\right)}{2}$ are raised to
the powers $1,2,3$ and we notice the simple pattern of 
\[
4=2^{2},\,8=2^{3},\,16=2^{4}.
\]
 On the left side we have sums involving $O_{3},\,O_{5},\,O_{7},$
and $O_{9}$ and coefficients in integers. We might not be able to
explain the integers yet, but we are confident enough to conjecture
the next expression will be 
\begin{equation}
o_{5}O_{5}+o_{7}O_{7}+o_{9}O_{9}+o_{11}O_{11}=32\cdot\left(\frac{n(n+1)}{2}\right)^{4}\label{eq:33}
\end{equation}
 for some integers $o_{5},o_{7},o_{9},o_{11}$.

\section{Back to the Hunt 2}

We are so excited by our new conjectures we hardly can wait to test
them out. Let us try the one for even powers: 
\begin{equation}
e_{4}E_{4}+e_{6}E_{6}+e_{8}E_{8}+e_{10}E_{10}=2^{4}\cdot3\cdot\left(\frac{n(n+1)}{2}\right)^{4}\label{eq:34}
\end{equation}
 for some integers $e_{4},e_{6},e_{8},e_{10}$. Can we use it to find
$E_{10}$?

\subsection{$E_{10}$}

From our previous work we suspect $E_{10}$ will have a form like
\begin{equation}
E_{10}=\frac{1}{11}\cdot\left(a_{1}\cdot\left(\frac{n(n+1)}{2}\right)^{4}-a_{2}\cdot\left(\frac{n(n+1)}{2}\right)^{3}+a_{3}\cdot\left(\frac{n(n+1)}{2}\right)^{2}-a_{4}\cdot E_{4}\right),\label{eq:35}
\end{equation}
 where $a_{1},a_{2},a_{3},a_{4}$ are rational numbers and $E_{4}$
is the familiar 
\[
E_{4}=\frac{1}{5}\cdot\left(6\cdot\frac{n(n+1)}{2}-1\right).
\]
 The coefficient of $2^{4}\cdot3$ in expression \ref{eq:34} tells
us $a_{1}=48$. Looking at expression \ref{eq:35} and the earlier
list in (\ref{eq:27}), we choose $e_{10}=11$. Therefore we rewrite
expression \ref{eq:34} as 
\begin{equation}
e_{4}E_{4}+e_{6}E_{6}+e_{8}E_{8}+11E_{10}=48\cdot\left(\frac{n(n+1)}{2}\right)^{4}\label{eq:36}
\end{equation}
 and expression \ref{eq:35} as 
\begin{equation}
E_{10}=\frac{1}{11}\cdot\left(48\cdot\left(\frac{n(n+1)}{2}\right)^{4}-a_{2}\cdot\left(\frac{n(n+1)}{2}\right)^{3}+a_{3}\cdot\left(\frac{n(n+1)}{2}\right)^{2}-a_{4}\cdot E_{4}\right).\label{eq:37}
\end{equation}

What does expression \ref{eq:36} tell us? On the right side of the
expression we see only $48\cdot\left(\frac{n(n+1)}{2}\right)^{4}$.
By expression \ref{eq:37} we \textit{suspect} $E_{10}$ has terms
of the form $\left(\frac{n(n+1)}{2}\right)^{4},\left(\frac{n(n+1)}{2}\right)^{3},\left(\frac{n(n+1)}{2}\right)^{2},$
and $E_{4}$. We \textit{know} $E_{4},E_{6},$ and $E_{8}$ have terms
of the form $\left(\frac{n(n+1)}{2}\right)^{3},\left(\frac{n(n+1)}{2}\right)^{2},$
and $E_{4}.$ That means all of the other terms must cancel out. In
other words, we have the following system of equations: 
\begin{align*}
E_{10}: &  & -a_{2}\cdot\left(\frac{n(n+1)}{2}\right)^{3} & +a_{3}\cdot\left(\frac{n(n+1)}{2}\right)^{2} & -a_{4}\cdot E_{4}\\
E_{8}: &  & e_{8}\cdot\frac{1}{9}\cdot24\cdot\left(\frac{n(n+1)}{2}\right)^{3} & -e_{8}\cdot\frac{1}{9}\cdot24\cdot\left(\frac{n(n+1)}{2}\right)^{2} & +e_{8}\cdot\frac{1}{9}\cdot9\cdot E_{4}\\
E_{6}: &  &  & e_{6}\cdot\frac{1}{7}\cdot12\cdot\left(\frac{n\left(n+1\right)}{2}\right)^{2} & -e_{6}\cdot\frac{1}{7}\cdot5\cdot E_{4}\\
E_{4}: &  &  &  & e_{4}\cdot E_{4}
\end{align*}
 where 
\begin{eqnarray*}
 & -a_{2}+e_{8}\cdot\frac{1}{9}\cdot24=0\\
 & a_{3}-e_{8}\cdot\frac{1}{9}\cdot24+e_{6}\cdot\frac{1}{7}\cdot12=0\\
 & -a_{4}+e_{8}\cdot\frac{1}{9}\cdot9-e_{6}\cdot\frac{1}{7}\cdot5+e_{4}=0,
\end{eqnarray*}
 which we simplify to 
\begin{eqnarray}
\frac{8}{3}\cdot e_{8} & = & a_{2}\nonumber \\
-\frac{12}{7}\cdot e_{6}+\frac{8}{3}\cdot e_{8} & = & a_{3}\nonumber \\
e_{4}-\frac{5}{7}\cdot e_{6}+e_{8} & = & a_{4}.\label{eq:38}
\end{eqnarray}
 How do we solve such a system of equations?

\bigskip{}
First, we need to decide whether we solve for $e_{i}$ or $a_{j}$.
That is easy to answer. If we knew the $a_{j}$ then we would have
solved the problem already. Second, that means if we substitute the
values for $e_{i}$ then we will get the values for $a_{j}$. What
\textit{are} the values for $e_{4},e_{6},$ and $e_{8}$?

If we look back to the previous expressions then we will remind ourselves
we \textit{derived }
\begin{eqnarray}
 & E_{2}=1\nonumber \\
 & E_{2}+5E_{4}=6\cdot\frac{n\left(n+1\right)}{2}\nonumber \\
 & 5E_{4}+7E_{6}=12\cdot\left(\frac{n\left(n+1\right)}{2}\right)^{2}\nonumber \\
 & E_{4}+14E_{6}+9E_{8}=24\cdot\left(\frac{n(n+1)}{2}\right)^{3}\label{eq:39}
\end{eqnarray}
 and \textit{guessed 
\[
e_{4}E_{4}+e_{6}E_{6}+e_{8}E_{8}+11E_{10}=48\cdot\left(\frac{n(n+1)}{2}\right)^{4}.
\]
} We \textit{do not know} the values for $e_{4},e_{6},$ and $e_{8}$.
What do we do?

Let us start by looking at the equations in (\ref{eq:38}). We observe
the following:
\begin{enumerate}
\item We have $\frac{8}{3}\cdot e_{8}=a_{2}$, which is equivalent to $e_{8}=\frac{3a_{2}}{8}$.
Since we suspect $e_{8}$ is an integer, if $a_{2}$ \textit{also}
is an integer then it follows that $a_{2}$ is a multiple of 8: $a_{2}=8c$
for some integer $c$. Therefore we may write $e_{8}=\frac{3\cdot8c}{8}=3c$,
which implies $e_{8}$ is a multiple of 3.\footnote{For another observation concerning $e_{8}$, see Section \ref{sub:Does--divide}
of the Appendix.}
\item We have $-\frac{12}{7}\cdot e_{6}+\frac{8}{3}\cdot e_{8}=a_{3}$,
which is equivalent to $e_{6}=\frac{7}{12}\cdot\left(\frac{8}{3}\cdot e_{8}-a_{3}\right)$.
Analogous to the case for $e_{8}$, if both $e_{6}$ and $a_{3}$
are integers then it follows that $e_{6}$ is a multiple of 7: $e_{6}=7d$
for some integer $d$.
\item We have $e_{4}-\frac{5}{7}\cdot e_{6}+e_{8}=a_{4}$, which is equivalent
to $e_{4}=\frac{5}{7}\cdot e_{6}-e_{8}+a_{4}$. That says little about
$e_{4}$. All it seems to suggest is that $a_{4}$ might be an integer
too.
\end{enumerate}
This is little to go on. What do we do?

Let us return to the expressions in (\ref{eq:39}), which suggested
the conjecture originally. However, like the approach of Section \ref{sec:A-Different-Point-of-View},
let us focus on the coefficients of the expressions rather than the
total expressions: 
\begin{align}
1 & =1\nonumber \\
6 & =1+5\nonumber \\
12 & =0+5+7\nonumber \\
24 & =0+1+14+9.\label{eq:40}
\end{align}
 Remarkable! In each expression the sum of the coefficients on the
right side is equal to the coefficient on the left side. Does it follow
that 
\begin{equation}
e_{4}+e_{6}+e_{8}+11=48,\label{eq:41}
\end{equation}
 which is equivalent to 
\begin{equation}
e_{4}+e_{6}+e_{8}=37\,?\label{eq:42}
\end{equation}
 If so, will it allow us to solve the previous system of equations?

Wait, there is more! Suppose we add the terms from expression \ref{eq:41}
to the list in (\ref{eq:40}). Then we get 
\begin{align}
1 & =1\nonumber \\
6 & =1+5\nonumber \\
12 & =0+5+7\nonumber \\
24 & =0+1+14+9\nonumber \\
48 & =0+e_{4}+e_{6}+e_{8}+11.\label{eq:43}
\end{align}
 Look at the terms on the diagonal containing $e_{8}$: 1,5,14. Have
we seen them before? Of course! They are 
\[
1=\sum_{k=1}^{1}k^{2},\,5=\sum_{k=1}^{2}k^{2},\,14=\sum_{k=1}^{3}k^{2}.
\]
 Does it follow that 
\[
e_{8}=30=\sum_{k=1}^{4}k^{2}\,?
\]
 Remembering our previous observations on the list of expressions
in (\ref{eq:38}), we noticed that $e_{8}$ is a multiple of 3. 30
is a multiple of 3! If we substitute $e_{8}=30$ into expression \ref{eq:42}
then we get $e_{4}+e_{6}=7$. We noticed also that $e_{6}$ is a multiple
of 7. 7 is a multiple of 7! Suppose we choose $e_{6}=7$, which implies
$e_{4}=0$. Then we may rewrite the list in (\ref{eq:43}) as 
\begin{align}
1 & =1\nonumber \\
6 & =1+5\nonumber \\
12 & =0+5+7\nonumber \\
24 & =0+1+14+9\nonumber \\
48 & =0+0+7+30+11\label{eq:44}
\end{align}
 and expression \ref{eq:36} as 
\begin{equation}
7E_{6}+30E_{8}+11E_{10}=48\cdot\left(\frac{n(n+1)}{2}\right)^{4}.\label{eq:45}
\end{equation}

Did we find the correct values for $e_{4},e_{6},$ and $e_{8}$? If
we substitute them into the equations of (\ref{eq:38}) then we get
\begin{align*}
a_{2} & =\frac{8}{3}\cdot30=80\\
a_{3} & =-\frac{12}{7}\cdot7+\frac{8}{3}\cdot30=68\\
a_{4} & =0-\frac{5}{7}\cdot7+30=25.
\end{align*}
 If we substitute these values for $a_{2},a_{3},$ and $a_{4}$ into
expression \ref{eq:37} then we get 
\begin{equation}
E_{10}=\frac{1}{11}\cdot\left(48\cdot\left(\frac{n(n+1)}{2}\right)^{4}-80\cdot\left(\frac{n(n+1)}{2}\right)^{3}+68\cdot\left(\frac{n(n+1)}{2}\right)^{2}-25\cdot E_{4}\right).\label{eq:46}
\end{equation}
 Last, if we return to the general form for sums of even powers then
we get 
\begin{equation}
11\cdot\frac{\sum k^{10}}{\sum k^{2}}=48\cdot\left(\frac{n(n+1)}{2}\right)^{4}-80\cdot\left(\frac{n(n+1)}{2}\right)^{3}+68\cdot\left(\frac{n(n+1)}{2}\right)^{2}-25\cdot\frac{1}{5}\cdot\left(6\cdot\frac{n\left(n+1\right)}{2}-1\right),\label{eq:47}
\end{equation}
 which actually gives us a way to test the new result: 
\[
\sum_{k=1}^{10}k^{10}=14,914,341,925
\]
 and 
\[
\frac{1}{11}\cdot\left(48\cdot55^{4}-80\cdot55^{3}+68\cdot55^{2}-25\cdot\frac{1}{5}\cdot\left(6\cdot55-1\right)\right)\cdot\sum_{k=1}^{10}k^{2}
\]
\[
=\frac{1}{11}\cdot426124055\cdot385=14,914,341,925.
\]
 It's only the single case of $n=10$, but that's good enough for
now.

\subsection{$O_{11}$}

Let us try the conjecture for the next odd power, $O_{11}$: 
\begin{equation}
o_{5}O_{5}+o_{7}O_{7}+o_{9}O_{9}+o_{11}O_{11}=2^{5}\cdot\left(\frac{n(n+1)}{2}\right)^{4}\label{eq:48}
\end{equation}
 for some integers $o_{5},o_{7},o_{9},o_{11}$, where 
\begin{equation}
o_{5}+o_{7}+o_{9}+o_{11}=2^{5}=32.\label{eq:49}
\end{equation}
 Since this case is a bit different than the previous one for $E_{10}$,
we remind ourselves that 
\[
O_{3}=1
\]
\[
O_{5}=\frac{1}{3}\cdot\left(4\cdot\frac{n\left(n+1\right)}{2}-1\right)
\]
\[
O_{7}=\frac{1}{3}\cdot\left(6\cdot\left(\frac{n(n+1)}{2}\right)^{2}-3\cdot O_{5}\right)
\]
\[
O_{9}=\frac{1}{10}\cdot\left(32\cdot\left(\frac{n(n+1)}{2}\right)^{3}-40\cdot\left(\frac{n(n+1)}{2}\right)^{2}+18\cdot O_{5}\right)
\]
 and 
\begin{eqnarray}
 & O_{3}=1\nonumber \\
 & O_{3}+3O_{5}=4\cdot\frac{n\left(n+1\right)}{2}\nonumber \\
 & 4O_{5}+4O_{7}=8\cdot\left(\frac{n(n+1)}{2}\right)^{2}\nonumber \\
 & O_{5}+10\cdot O_{7}+5O_{9}=16\cdot\left(\frac{n(n+1)}{2}\right)^{3}.\label{eq:50}
\end{eqnarray}

We suspect $O_{11}$ will have a form like 
\begin{equation}
O_{11}=\frac{1}{6}\cdot\left(32\cdot\left(\frac{n(n+1)}{2}\right)^{4}-b_{2}\cdot\left(\frac{n(n+1)}{2}\right)^{3}+b_{3}\cdot\left(\frac{n(n+1)}{2}\right)^{2}-b_{4}\cdot O_{5}\right)\label{eq:51}
\end{equation}
 for some rational numbers $b_{2},b_{3},b_{4}$. Therefore we rewrite
expression \ref{eq:48} as 
\begin{equation}
o_{5}O_{5}+o_{7}O_{7}+o_{9}O_{9}+6\cdot O_{11}=32\cdot\left(\frac{n(n+1)}{2}\right)^{4}\label{eq:52}
\end{equation}
 and expression \ref{eq:49} as 
\begin{equation}
o_{5}+o_{7}+o_{9}+6=32,\label{eq:53}
\end{equation}
 which is equivalent to 
\begin{equation}
o_{5}+o_{7}+o_{9}=26.\label{eq:54}
\end{equation}
 Let us gather more information about $o_{5},o_{7},$ and $o_{9}$.

\bigskip{}
From the list in (\ref{eq:50}) we know that 
\begin{align*}
1 & =1\\
4 & =1+3\\
8 & =0+4+4\\
16 & =0+1+10+5.
\end{align*}
 If we add expression \ref{eq:53} to it the we get 
\begin{align}
1 & =1\nonumber \\
4 & =1+3\nonumber \\
8 & =0+4+4\nonumber \\
16 & =0+1+10+5\nonumber \\
32 & =0+o_{5}+o_{7}+o_{9}+6.\label{eq:55}
\end{align}
 What are the terms along the diagonal for $o_{9}$? We suspect 
\[
4-1=3=1+2=\sum_{k=1}^{2}k
\]
\[
10-4=6=1+2+3=\sum_{k=1}^{3}k
\]
 and 
\[
o_{9}-10=10=1+2+3+4=\sum_{k=1}^{4}k,
\]
 which implies $o_{9}=20$. Therefore we rewrite the list in (\ref{eq:55})
as 
\begin{align}
1 & =1\nonumber \\
4 & =1+3\nonumber \\
8 & =0+4+4\nonumber \\
16 & =0+1+10+5\nonumber \\
32 & =0+o_{5}+o_{7}+20+6\label{eq:56}
\end{align}
 and expression \ref{eq:52} as 
\[
o_{5}O_{5}+o_{7}O_{7}+20\cdot O_{9}+6\cdot O_{11}=32\cdot\left(\frac{n(n+1)}{2}\right)^{4}.
\]
 That still leaves 
\[
o_{5}+o_{7}=6
\]
 and 
\[
O_{11}=\frac{1}{6}\cdot\left(32\cdot\left(\frac{n(n+1)}{2}\right)^{4}-b_{2}\cdot\left(\frac{n(n+1)}{2}\right)^{3}+b_{3}\cdot\left(\frac{n(n+1)}{2}\right)^{2}-b_{4}\cdot O_{5}\right).
\]
 How do we figure out $o_{5},o_{7}$ and $b_{2},b_{3},b_{4}$?

\bigskip{}
In order to determine all of the values for $o_{i}$ and $b_{j}$
we need to solve the following system of equations: 
\begin{align*}
O_{11}: &  & -b_{2}\cdot\left(\frac{n(n+1)}{2}\right)^{3} & +b_{3}\cdot\left(\frac{n(n+1)}{2}\right)^{2} & -b_{4}\cdot O_{5}\\
O_{9}: &  & o_{9}\cdot\frac{32}{10}\cdot\left(\frac{n(n+1)}{2}\right)^{3} & -o_{9}\cdot\frac{40}{10}\cdot\left(\frac{n(n+1)}{2}\right)^{2} & +o_{9}\cdot\frac{18}{10}\cdot O_{5}\\
O_{7}: &  &  & o_{7}\cdot\frac{6}{3}\cdot\left(\frac{n(n+1)}{2}\right)^{2} & -o_{7}\cdot\frac{3}{3}\cdot O_{5}\\
O_{5}: &  &  &  & o_{5}\cdot O_{5},
\end{align*}
 where 
\begin{eqnarray*}
 & -b_{2}+o_{9}\cdot\frac{32}{10}=0\\
 & b_{3}-o_{9}\cdot\frac{40}{10}+o_{7}\cdot\frac{6}{3}=0\\
 & -b_{4}+o_{9}\cdot\frac{18}{10}-o_{7}\cdot\frac{3}{3}+o_{5}=0,
\end{eqnarray*}
 which we simplify to 
\begin{eqnarray}
o_{9}\cdot\frac{16}{5} & = & b_{2}\nonumber \\
o_{9}\cdot4-o_{7}\cdot2 & = & b_{3}\nonumber \\
o_{9}\cdot\frac{9}{5}-o_{7}+o_{5} & = & b_{4}.\label{eq:57}
\end{eqnarray}
 $o_{9}=20$ implies 
\[
b_{2}=20\cdot\frac{16}{5}=64,
\]
 which allows us to rewrite expression \ref{eq:51} as 
\begin{equation}
O_{11}=\frac{1}{6}\cdot\left(32\cdot\left(\frac{n(n+1)}{2}\right)^{4}-64\cdot\left(\frac{n(n+1)}{2}\right)^{3}+b_{3}\cdot\left(\frac{n(n+1)}{2}\right)^{2}-b_{4}\cdot O_{5}\right).\label{eq:58}
\end{equation}
 As for the remaining variables, from the equations in (\ref{eq:57})
we observe
\begin{enumerate}
\item $o_{9}\cdot\frac{16}{5}=b_{2}$, which is equivalent to $o_{9}=\frac{5b_{2}}{16}$.
If both $o_{9}$ and $b_{2}$ are integers then $b_{2}$ is a multiple
of 16, which implies $o_{9}$ is a multiple of 5: $o_{9}=5c$ for
some integer $c$. 
\item $o_{9}\cdot4-o_{7}\cdot2=b_{3}$, which is equivalent to $2\cdot\left(o_{7}-2\cdot o_{9}\right)=b_{3}$.
If both $o_{7}$ and $o_{9}$ are integers then $b_{3}$ is even.
\item $o_{9}\cdot\frac{9}{5}-o_{7}+o_{5}=b_{4}$ might imply $b_{4}$ is
an integer. 
\end{enumerate}
Once again, that is little to go on. All it suggests is the guess
of $o_{9}=20$ might be correct. What do we do?

Suppose we make the simple choice of $o_{7}=6$, which implies $o_{5}=0$.
Then we may rewrite the list in (\ref{eq:56}) as 
\begin{align}
1 & =1\nonumber \\
4 & =1+3\nonumber \\
8 & =0+4+4\nonumber \\
16 & =0+1+10+5\nonumber \\
32 & =0+0+6+20+6\label{eq:59}
\end{align}
 and expression \ref{eq:52} as 
\begin{equation}
6\cdot O_{7}+20\cdot O_{9}+6\cdot O_{11}=32\cdot\left(\frac{n(n+1)}{2}\right)^{4}.\label{eq:60}
\end{equation}

Did we find the correct values for $o_{5},o_{7}$ and $o_{9}$? If
we substitute them into the equations of (\ref{eq:57}) then we get
\begin{align*}
b_{3} & =20\cdot4-6\cdot2=68\\
b_{4} & =20\cdot\frac{9}{5}-6+0=30.
\end{align*}
 If we substitute these values for $b_{3}$ and $b_{4}$ into expression
\ref{eq:58} then we get 
\begin{equation}
O_{11}=\frac{1}{6}\cdot\left(32\cdot\left(\frac{n(n+1)}{2}\right)^{4}-64\cdot\left(\frac{n(n+1)}{2}\right)^{3}+68\cdot\left(\frac{n(n+1)}{2}\right)^{2}-30\cdot O_{5}\right).\label{eq:61}
\end{equation}
 Last, if we return to the general form for sums of odd powers then
we get 
\begin{equation}
6\cdot\frac{\sum k^{11}}{\sum k^{3}}=32\cdot\left(\frac{n(n+1)}{2}\right)^{4}-64\cdot\left(\frac{n(n+1)}{2}\right)^{3}+68\cdot\left(\frac{n(n+1)}{2}\right)^{2}-30\cdot\frac{1}{3}\cdot\left(4\cdot\frac{n\left(n+1\right)}{2}-1\right),\label{eq:62}
\end{equation}
 which actually gives us a way to test the new result: 
\[
\sum_{k=1}^{9}k^{11}=42,364,319,625
\]
 and 
\[
\frac{1}{6}\cdot\left(32\cdot45^{4}-64\cdot45^{3}+68\cdot45^{2}-30\cdot\frac{1}{3}\cdot\left(4\cdot45-1\right)\right)\cdot\sum_{k=1}^{9}k^{3}
\]
 
\[
=\frac{1}{6}\cdot125523910\cdot2025=42,364,319,625.
\]
 We'll take it.

\section{A Chance to Catch Our Breath}

Our new results for even powers, which are only conjectures, are 
\[
7E_{6}+30E_{8}+11E_{10}=48\cdot\left(\frac{n(n+1)}{2}\right)^{4}
\]
\[
E_{10}=\frac{1}{11}\cdot\left(48\cdot\left(\frac{n(n+1)}{2}\right)^{4}-80\cdot\left(\frac{n(n+1)}{2}\right)^{3}+68\cdot\left(\frac{n(n+1)}{2}\right)^{2}-25\cdot E_{4}\right)
\]
\[
11\cdot\frac{\sum k^{10}}{\sum k^{2}}=48\cdot\left(\frac{n(n+1)}{2}\right)^{4}-80\cdot\left(\frac{n(n+1)}{2}\right)^{3}+68\cdot\left(\frac{n(n+1)}{2}\right)^{2}-25\cdot\frac{1}{5}\cdot\left(6\cdot\frac{n\left(n+1\right)}{2}-1\right).
\]
 They allow us to expand our list to 
\begin{eqnarray*}
 & E_{2}=1\\
 & E_{2}+5E_{4}=2\cdot3\cdot\frac{n\left(n+1\right)}{2}\\
 & 5E_{4}+7E_{6}=2^{2}\cdot3\cdot\left(\frac{n\left(n+1\right)}{2}\right)^{2}\\
 & E_{4}+14E_{6}+9E_{8}=2^{3}\cdot3\cdot\left(\frac{n(n+1)}{2}\right)^{3}\\
 & 7E_{6}+30E_{8}+11E_{10}=2^{4}\cdot3\cdot\left(\frac{n(n+1)}{2}\right)^{4},
\end{eqnarray*}
 which we find more advantageous to express as 
\begin{align*}
1 & =1\\
6 & =1+5\\
12 & =0+5+7\\
24 & =0+1+14+9\\
48 & =0+0+7+30+11.
\end{align*}

Our new results for odd powers, which also are only conjectures, are
\[
6\cdot O_{7}+20\cdot O_{9}+6\cdot O_{11}=32\cdot\left(\frac{n(n+1)}{2}\right)^{4}
\]
\[
O_{11}=\frac{1}{6}\cdot\left(32\cdot\left(\frac{n(n+1)}{2}\right)^{4}-64\cdot\left(\frac{n(n+1)}{2}\right)^{3}+68\cdot\left(\frac{n(n+1)}{2}\right)^{2}-30\cdot O_{5}\right)
\]
\[
6\cdot\frac{\sum k^{11}}{\sum k^{3}}=32\cdot\left(\frac{n(n+1)}{2}\right)^{4}-64\cdot\left(\frac{n(n+1)}{2}\right)^{3}+68\cdot\left(\frac{n(n+1)}{2}\right)^{2}-30\cdot\frac{1}{3}\cdot\left(4\cdot\frac{n\left(n+1\right)}{2}-1\right).
\]
 They allow us to expand our list to 
\begin{eqnarray*}
 & O_{3}=1\\
 & O_{3}+3O_{5}=2^{2}\cdot\frac{n\left(n+1\right)}{2}\\
 & 4O_{5}+4O_{7}=2^{3}\cdot\left(\frac{n(n+1)}{2}\right)^{2}\\
 & O_{5}+10\cdot O_{7}+5O_{9}=2^{4}\cdot\left(\frac{n(n+1)}{2}\right)^{3}\\
 & 6\cdot O_{7}+20\cdot O_{9}+6\cdot O_{11}=2^{5}\cdot\left(\frac{n(n+1)}{2}\right)^{4},
\end{eqnarray*}
 which we find more advantageous to express as 
\begin{align*}
1 & =1\\
4 & =1+3\\
8 & =0+4+4\\
16 & =0+1+10+5\\
32 & =0+0+6+20+6.
\end{align*}
 The outer terms suggest we recast some previous expressions for $O_{2m+1}$:
\begin{eqnarray*}
O_{7} & = & \frac{1}{3}\cdot\left(6\cdot\left(\frac{n(n+1)}{2}\right)^{2}-3\cdot O_{5}\right)\\
 &  & =\frac{1}{4}\cdot\left(8\cdot\left(\frac{n(n+1)}{2}\right)^{2}-4\cdot O_{5}\right)
\end{eqnarray*}
\begin{eqnarray*}
O_{9} & = & \frac{1}{10}\cdot\left(32\cdot\left(\frac{n(n+1)}{2}\right)^{3}-40\cdot\left(\frac{n(n+1)}{2}\right)^{2}+18\cdot O_{5}\right)\\
 &  & =\frac{1}{5}\cdot\left(16\cdot\left(\frac{n(n+1)}{2}\right)^{3}-20\cdot\left(\frac{n(n+1)}{2}\right)^{2}+9\cdot O_{5}\right).
\end{eqnarray*}
\pagebreak{}

\section{Summary of Part 3}

The result of our hard work in Parts 1 and 2 were explicit expressions
for 
\[
\sum k^{2},\,\sum k^{3},\ldots,\sum k^{8},\,\sum k^{9},
\]
 which we divided into even and odd powers. Unfortunately, after inspecting
the list of sums we were unable to discover any general patterns.
Following up on an earlier observation, we came across the idea that
in expressing the sums as 
\[
\sum k^{2m}=E_{2m}\cdot\sum k^{2}
\]
\[
\sum k^{2m+1}=O_{2m+1}\cdot\sum k^{3},
\]
 perhaps some insight was to be gained by removing the coefficients
$E_{2m}$ and $O_{2m+1}$ and looking at their relationships amongst
only themselves.

To our shock we discovered the equations 
\begin{eqnarray*}
 & E_{2}=1\\
 & E_{2}+5E_{4}=2\cdot3\cdot\frac{n\left(n+1\right)}{2}\\
 & 5E_{4}+7E_{6}=2^{2}\cdot3\cdot\left(\frac{n\left(n+1\right)}{2}\right)^{2}\\
 & E_{4}+14E_{6}+9E_{8}=2^{3}\cdot3\cdot\left(\frac{n(n+1)}{2}\right)^{3}
\end{eqnarray*}
 and 
\begin{eqnarray*}
 & O_{3}=1\\
 & O_{3}+3O_{5}=2^{2}\cdot\frac{n\left(n+1\right)}{2}\\
 & 4O_{5}+4O_{7}=2^{3}\cdot\left(\frac{n(n+1)}{2}\right)^{2}\\
 & O_{5}+10\cdot O_{7}+5O_{9}=2^{4}\cdot\left(\frac{n(n+1)}{2}\right)^{3},
\end{eqnarray*}
 which led us to conjecture 
\[
e_{4}E_{4}+e_{6}E_{6}+e_{8}E_{8}+e_{10}E_{10}=2^{4}\cdot3\cdot\left(\frac{n(n+1)}{2}\right)^{4}
\]
\[
o_{5}O_{5}+o_{7}O_{7}+o_{9}O_{9}+o_{11}O_{11}=2^{5}\cdot\left(\frac{n(n+1)}{2}\right)^{4}.
\]
 Through inductive reasoning and a simple system of equations we derived
\[
7E_{6}+30E_{8}+11E_{10}=48\cdot\left(\frac{n(n+1)}{2}\right)^{4}
\]
\[
6\cdot O_{7}+20\cdot O_{9}+6\cdot O_{11}=32\cdot\left(\frac{n(n+1)}{2}\right)^{4},
\]
 which yielded 
\[
E_{10}=\frac{1}{11}\cdot\left(48\cdot\left(\frac{n(n+1)}{2}\right)^{4}-80\cdot\left(\frac{n(n+1)}{2}\right)^{3}+68\cdot\left(\frac{n(n+1)}{2}\right)^{2}-25\cdot E_{4}\right)
\]
\[
O_{11}=\frac{1}{6}\cdot\left(32\cdot\left(\frac{n(n+1)}{2}\right)^{4}-64\cdot\left(\frac{n(n+1)}{2}\right)^{3}+68\cdot\left(\frac{n(n+1)}{2}\right)^{2}-30\cdot O_{5}\right).
\]
 Even though the new results were only conjectures, by reinserting
them into 
\[
\sum k^{10}=E_{10}\cdot\sum k^{2}
\]
\[
\sum k^{11}=O_{11}\cdot\sum k^{3}
\]
 we were able to verify them in special cases. Given how complicated
the sums were, we had good reason to believe we found the correct
expressions.\footnote{For an example of when inductive reasoning goes wrong, see Section
\ref{sub:Inductive-Reasoning-Sometimes} of the Appendix.}

Last, analogous to how in the expressions for the sums the coefficients
turned out to hold the more important relationships, we believe that
in the expressions for the coefficients \textit{their} coefficients
will turn out to hold the more important relationships, namely 
\begin{align*}
1 & =1\\
6 & =1+5\\
12 & =0+5+7\\
24 & =0+1+14+9\\
48 & =0+0+7+30+11
\end{align*}
 and 
\begin{align*}
1 & =1\\
4 & =1+3\\
8 & =0+4+4\\
16 & =0+1+10+5\\
32 & =0+0+6+20+6.
\end{align*}

\pagebreak{}

\part*{Part 4}

\section{A Frenchman Comes Aboard\label{sec:A-Frenchman-Comes}}

At the close of Part 3 we placed our hopes in new observations about
$E_{2m}$ and $O_{2m+1}$, the expressions 
\begin{eqnarray*}
 & E_{2}=1\\
 & E_{2}+5E_{4}=2\cdot3\cdot\frac{n\left(n+1\right)}{2}\\
 & 5E_{4}+7E_{6}=2^{2}\cdot3\cdot\left(\frac{n\left(n+1\right)}{2}\right)^{2}\\
 & E_{4}+14E_{6}+9E_{8}=2^{3}\cdot3\cdot\left(\frac{n(n+1)}{2}\right)^{3}\\
 & 7E_{6}+30E_{8}+11E_{10}=2^{4}\cdot3\cdot\left(\frac{n(n+1)}{2}\right)^{4}
\end{eqnarray*}
 and 
\begin{eqnarray*}
 & O_{3}=1\\
 & O_{3}+3O_{5}=2^{2}\cdot\frac{n\left(n+1\right)}{2}\\
 & 4O_{5}+4O_{7}=2^{3}\cdot\left(\frac{n(n+1)}{2}\right)^{2}\\
 & O_{5}+10\cdot O_{7}+5O_{9}=2^{4}\cdot\left(\frac{n(n+1)}{2}\right)^{3}\\
 & 6O_{7}+20\cdot O_{9}+6O_{11}=2^{5}\cdot\left(\frac{n(n+1)}{2}\right)^{4},
\end{eqnarray*}
 and the relationships amongst their coefficients, 
\begin{align}
1 & =1\nonumber \\
2\cdot3 & =1+5\nonumber \\
2^{2}\cdot3 & =0+5+7\nonumber \\
2^{3}\cdot3 & =0+1+14+9\nonumber \\
2^{4}\cdot3 & =0+0+7+30+11\label{eq:63}
\end{align}
 and 
\begin{align}
1 & =1\nonumber \\
2^{2} & =1+3\nonumber \\
2^{3} & =0+4+4\nonumber \\
2^{4} & =0+1+10+5\nonumber \\
2^{5} & =0+0+6+20+6.\label{eq:64}
\end{align}
 We believed they held the key to discovering general expressions
for $E_{2m}$ and $O_{2m+1}$. To continue the investigation, from
the outward appearances of (\ref{eq:63}) and (\ref{eq:64}) we hardly
can do any better than to look for relationships in Pascal's Triangle.

\bigskip{}
To construct Pascal's Triangle we start with a 1 at the top and a
1,2,1 in the second row, place a 1 at the start and end of each subsequent
row, and fill in the middle of the triangle by adding adjacent terms:
\[
3=1+2,\,3=2+1
\]
\[
4=1+3,\,6=3+3,\,4=3+1
\]
\[
5=1+4,\,10=4+6,\,10=6+4,\,5=4+1,
\]
 and so forth. We get

\begin{longtable}[l]{ccccccccccccccccc}
 &  &  &  &  &  &  &  & 1 &  &  &  &  &  &  &  & \tabularnewline
 &  &  &  &  &  &  &  &  &  &  &  &  &  &  &  & \tabularnewline
 &  &  &  &  &  &  &  &  &  &  &  &  &  &  &  & \tabularnewline
 &  &  &  &  &  & 1 &  & 2 &  & 1 &  &  &  &  &  & \tabularnewline
 &  &  &  &  &  &  &  &  &  &  &  &  &  &  &  & \tabularnewline
 &  &  &  &  &  &  &  &  &  &  &  &  &  &  &  & \tabularnewline
 &  &  &  &  & 1 &  & 3 &  & 3 &  & 1 &  &  &  &  & \tabularnewline
 &  &  &  &  &  &  &  &  &  &  &  &  &  &  &  & \tabularnewline
 &  &  &  &  &  &  &  &  &  &  &  &  &  &  &  & \tabularnewline
 &  &  &  & 1 &  & 4 &  & 6 &  & 4 &  & 1 &  &  &  & \tabularnewline
 &  &  &  &  &  &  &  &  &  &  &  &  &  &  &  & \tabularnewline
 &  &  &  &  &  &  &  &  &  &  &  &  &  &  &  & \tabularnewline
 &  &  & 1 &  & 5 &  & 10 &  & 10 &  & 5 &  & 1 &  &  & \tabularnewline
 &  &  &  &  &  &  &  &  &  &  &  &  &  &  &  & \tabularnewline
 &  &  &  &  &  &  &  &  &  &  &  &  &  &  &  & \tabularnewline
 &  & 1 &  & 6 &  & 15 &  & 20 &  & 15 &  & 6 &  & 1 &  & \tabularnewline
 &  &  &  &  &  &  &  &  &  &  &  &  &  &  &  & \tabularnewline
 &  &  &  &  &  &  &  &  &  &  &  &  &  &  &  & \tabularnewline
 & 1 &  & 7 &  & 21 &  & 35 &  & 35 &  & 21 &  & 7 &  & 1 & \tabularnewline
 &  &  &  &  &  &  &  &  &  &  &  &  &  &  &  & \tabularnewline
 &  &  &  &  &  &  &  &  &  &  &  &  &  &  &  & \tabularnewline
1 &  & 8 &  & 28 &  & 56 &  & 70 &  & 56 &  & 28 &  & 8 &  & 1\tabularnewline
\end{longtable} We may continue the process for as long as we like. To refer to entries
in the triangle we use the the notation $\binom{n}{k}$ . For example,
\[
\binom{1}{1}=1,\,\binom{4}{0}=1,\,\binom{4}{1}=2,\,\binom{4}{4}=4,\,\binom{7}{2}=21.
\]
 Entries outside the table we set equal to zero: for example, $\binom{6}{-1}=0$
and $\binom{2}{10}=0$. Does Pascal's Triangle help us to notice any
new relationships in the lists in (\ref{eq:63}) and (\ref{eq:64})?

\bigskip{}
DOES IT EVER!? If we look at the coefficients for $O_{2m+1}$, it's
as if the entire solution has been laid bare before our eyes:

\begin{align*}
1 & =1\\
2^{2} & =\binom{3}{0}+\binom{3}{2}\\
2^{3} & =\binom{4}{1}+\binom{4}{3}\\
2^{4} & =\binom{5}{0}+\binom{5}{2}+\binom{5}{4}\\
2^{5} & =\binom{6}{1}+\binom{6}{3}+\binom{6}{5}.
\end{align*}
 More important, from these few cases we believe we can extend the
list indefinitely: 
\begin{align*}
1 & =1\\
2^{2} & =\binom{3}{0}+\binom{3}{2}\\
2^{3} & =\binom{4}{1}+\binom{4}{3}\\
2^{4} & =\binom{5}{0}+\binom{5}{2}+\binom{5}{4}\\
2^{5} & =\binom{6}{1}+\binom{6}{3}+\binom{6}{5}\\
2^{6} & =\binom{7}{0}+\binom{7}{2}+\binom{7}{4}+\binom{7}{6}\\
2^{7} & =\binom{8}{1}+\binom{8}{3}+\binom{8}{5}+\binom{8}{7}\\
\vdots\\
2^{m} & =\binom{m+1}{0}+\binom{m+1}{2}+\binom{m+1}{4}+\cdots+\binom{m+1}{m}\\
2^{m+1} & =\binom{m+2}{1}+\binom{m+2}{3}+\binom{m+2}{5}+\cdots+\binom{m+2}{m+1}.
\end{align*}
 Of course sometimes it will be helpful to write it as 
\begin{align}
1 & =1\nonumber \\
4 & =1+3\nonumber \\
8 & =0+4+4\nonumber \\
16 & =0+1+10+5\nonumber \\
32 & =0+0+6+20+6\nonumber \\
64 & =0+0+1+21+35+7\nonumber \\
128 & =0+0+0+8+56+56+8.\label{eq:65}
\end{align}
 What about the coefficients for even powers, $E_{2m}$? The entries
are a bit harder to find, but we believe they are 
\begin{align*}
1 & =1\\
2\cdot3 & =\left[\binom{2}{-1}+\binom{3}{0}\right]+\left[\binom{2}{1}+\binom{3}{2}\right]\\
2^{2}\cdot3 & =\left[\binom{3}{0}+\binom{4}{1}\right]+\left[\binom{3}{2}+\binom{4}{3}\right]\\
2^{3}\cdot3 & =\left[\binom{4}{-1}+\binom{5}{0}\right]+\left[\binom{4}{1}+\binom{5}{2}\right]+\left[\binom{4}{3}+\binom{5}{4}\right]\\
2^{4}\cdot3 & =\left[\binom{5}{0}+\binom{6}{1}\right]+\left[\binom{5}{2}+\binom{6}{3}\right]+\left[\binom{5}{4}+\binom{6}{5}\right]\\
2^{5}\cdot3 & =\left[\binom{6}{-1}+\binom{7}{0}\right]+\left[\binom{6}{1}+\binom{7}{2}\right]+\left[\binom{6}{3}+\binom{7}{4}\right]+\left[\binom{6}{5}+\binom{7}{6}\right]\\
2^{6}\cdot3 & =\left[\binom{7}{0}+\binom{8}{1}\right]+\left[\binom{7}{2}+\binom{8}{3}\right]+\left[\binom{7}{4}+\binom{8}{5}\right]+\left[\binom{7}{6}+\binom{8}{7}\right]\\
\vdots\\
2^{m-1}\cdot3 & =\left[\binom{m}{-1}+\binom{m+1}{0}\right]+\left[\binom{m}{1}+\binom{m+1}{2}\right]+\left[\binom{m}{3}+\binom{m+1}{4}\right]\\
 & +\cdots+\left[\binom{m}{m-1}+\binom{m+1}{m}\right]\\
2^{m}\cdot3 & =\left[\binom{m+1}{0}+\binom{m+2}{1}\right]+\left[\binom{m+1}{2}+\binom{m+2}{3}\right]+\left[\binom{m+1}{4}+\binom{m+2}{5}\right]\\
 & +\cdots+\left[\binom{m+1}{m}+\binom{m+2}{m+1}\right],
\end{align*}
 which we still can write as 
\begin{align}
1 & =1\nonumber \\
6 & =1+5\nonumber \\
12 & =0+5+7\nonumber \\
24 & =0+1+14+9\nonumber \\
48 & =0+0+7+30+11\nonumber \\
96 & =0+0+1+27+55+13\nonumber \\
192 & =0+0+0+9+77+91+15.\label{eq:66}
\end{align}
 How did we deserve such luck?

\section{Disbelief}

Before we continue, let us look at the new discovery in more detail.
We observe the following:
\begin{enumerate}
\item put simply, we are at a loss for words about how the introduction
of Pascal's Triangle eliminated our troubles in one fell swoop. Is
there a deeper connection with the problem? Is there something we
missed? Also, we wonder why it took so long to make the connection.
Next time, in similar circumstances, will we be able to make such
a connection quicker?
\item when recast in terms of entries in Pascal's Triangle, the previous
relationships, which are rather ordinary, suggest new, concrete expressions.
For example, the sum 
\[
32=0+0+6+20+6
\]
 becomes 
\[
2^{5}=\binom{6}{1}+\binom{6}{3}+\binom{6}{5}.
\]
 The sum 
\[
48=0+0+7+30+11
\]
 becomes 
\[
2^{4}\cdot3=\left[\binom{5}{0}+\binom{6}{1}\right]+\left[\binom{5}{2}+\binom{6}{3}\right]+\left[\binom{5}{4}+\binom{6}{5}\right].
\]
Or, the sums $\sum k$ and $\sum k^{2}$ can be expressed as 
\begin{eqnarray*}
\sum_{k=1}^{5}k & = & 1+2+3+4+5=15\\
 &  & =1+\binom{2}{1}+\binom{3}{2}+\binom{4}{3}+\binom{5}{4}=\binom{6}{4};
\end{eqnarray*}
\begin{eqnarray*}
\sum_{k=1}^{5}k^{2} & = & 1+4+9+16+25=55\\
 &  & =20+35=\binom{6}{3}+\binom{7}{4}=\binom{6}{3}+\binom{7}{3};
\end{eqnarray*}
\begin{eqnarray*}
\sum_{k=1}^{1}k+\sum_{k=1}^{2}k+\sum_{k=1}^{3}k+\sum_{k=1}^{4}k+\sum_{k=1}^{5}k & = & 1+3+6+10+15=35\\
 &  & =1+\binom{3}{1}+\binom{4}{2}+\binom{5}{3}+\binom{6}{4}=\binom{7}{4};
\end{eqnarray*}
 
\[
\sum_{k=1}^{1}k^{2}+\sum_{k=1}^{2}k^{2}+\sum_{k=1}^{3}k^{2}+\sum_{k=1}^{4}k^{2}+\sum_{k=1}^{5}k^{2}=1+5+14+30+55=105
\]
\begin{eqnarray*}
 & = & 1+\left[\binom{3}{0}+\binom{4}{1}\right]+\left[\binom{4}{1}+\binom{5}{2}\right]+\left[\binom{5}{2}+\binom{6}{3}\right]+\left[\binom{6}{3}+\binom{7}{4}\right]\\
 &  & =35+70=\binom{7}{3}+\binom{8}{4}.
\end{eqnarray*}
 We believe we only have begun to scratch the surface. Unfortunately,
we must save such investigations for another time.
\end{enumerate}

\section{Back to the Hunt 3 }

Our new results lead to a substantial improvement in our ability to
calculate the desired expressions. Consider $E_{12}$ and $O_{13}.$

\subsection{$E_{12}$\label{sub:E12}}

From the list in (\ref{eq:66}) we conjecture 
\begin{equation}
E_{6}+27E_{8}+55E_{10}+13E_{12}=96\cdot\left(\frac{n(n+1)}{2}\right)^{5},
\end{equation}
 where 
\begin{eqnarray}
96 & = & e_{6}+e_{8}+e_{10}+e_{12}\nonumber \\
 &  & =1+27+55+13.
\end{eqnarray}

\pagebreak{}

We suspect $E_{12}$ will have the form 
\begin{equation}
13\cdot E_{12}=96\cdot\left(\frac{n(n+1)}{2}\right)^{5}-a_{2}\cdot\left(\frac{n(n+1)}{2}\right)^{4}+a_{3}\cdot\left(\frac{n(n+1)}{2}\right)^{3}-a_{4}\cdot\left(\frac{n(n+1)}{2}\right)^{2}+a_{5}\cdot E_{4}.\label{eq:69}
\end{equation}
 The system of equations is 
\begin{align*}
E_{12}: &  & -a_{2}\cdot\left(\frac{n(n+1)}{2}\right)^{4} & +a_{3}\cdot\left(\frac{n(n+1)}{2}\right)^{3} & -a_{4}\cdot\left(\frac{n(n+1)}{2}\right)^{2} & +a_{5}\cdot E_{4}\\
E_{10}: &  & 55\cdot\frac{48}{11}\cdot\left(\frac{n(n+1)}{2}\right)^{4} & -55\cdot\frac{80}{11}\cdot\left(\frac{n(n+1)}{2}\right)^{3} & +55\cdot\frac{68}{11}\cdot\left(\frac{n(n+1)}{2}\right)^{2} & -55\cdot\frac{25}{11}\cdot E_{4}\\
E_{8}: &  &  & 27\cdot\frac{24}{9}\cdot\left(\frac{n(n+1)}{2}\right)^{3} & -27\cdot\frac{24}{9}\cdot\left(\frac{n(n+1)}{2}\right)^{2} & +27\cdot\frac{9}{9}\cdot E_{4}\\
E_{6}: &  &  &  & 1\cdot\frac{12}{7}\cdot\left(\frac{n(n+1)}{2}\right)^{2} & -1\cdot\frac{5}{7}\cdot E_{4},
\end{align*}
 where 
\begin{eqnarray*}
 & -a_{2}+55\cdot\frac{48}{11}=0\\
 & a_{3}-55\cdot\frac{80}{11}+27\cdot\frac{24}{9}=0\\
 & -a_{4}+55\cdot\frac{68}{11}-27\cdot\frac{24}{9}+1\cdot\frac{12}{7}=0\\
 & a_{5}-55\cdot\frac{25}{11}+27\cdot\frac{9}{9}-1\cdot\frac{5}{7}=0.
\end{eqnarray*}
 We can solve it immediately: 
\begin{align*}
a_{2} & =240\\
a_{3} & =328\\
a_{4} & =\frac{1888}{7}\\
a_{5} & =\frac{691}{7}.
\end{align*}
 Therefore we can rewrite expression \ref{eq:69} as 
\begin{equation}
13\cdot E_{12}=96\cdot\left(\frac{n(n+1)}{2}\right)^{5}-240\cdot\left(\frac{n(n+1)}{2}\right)^{4}+328\cdot\left(\frac{n(n+1)}{2}\right)^{3}-\frac{1888}{7}\cdot\left(\frac{n(n+1)}{2}\right)^{2}+\frac{691}{7}\cdot E_{4}.
\end{equation}
 If we insert it into 
\[
\sum k^{12}=E_{12}\cdot\sum k^{2}
\]
 then we can test the result: 
\[
\sum_{k=1}^{9}k^{12}=367,428,536,133
\]
 and 
\[
\frac{1}{13}\cdot\left(96\cdot45^{5}-240\cdot45^{4}+328\cdot45^{3}-\frac{1888}{7}\cdot45^{2}+\frac{691}{7}\cdot\frac{1}{5}\cdot\left(6\cdot45-1\right)\right)\cdot\sum_{k=1}^{9}k^{2}
\]
\[
=\frac{1}{13}\cdot\frac{83799490697}{5}\cdot285=367,428,536,133.
\]

\subsection{$O_{13}$}

From the list in (\ref{eq:65}) we conjecture 
\begin{equation}
O_{7}+21\cdot O_{9}+35\cdot O_{11}+7\cdot O_{13}=64\cdot\left(\frac{n(n+1)}{2}\right)^{5},
\end{equation}
 where 
\begin{eqnarray}
64 & = & o_{7}+o_{9}+o_{11}+o_{13}\nonumber \\
 &  & =1+21+35+7.
\end{eqnarray}
 We suspect $O_{13}$ will have the form 
\begin{equation}
7\cdot O_{13}=64\cdot\left(\frac{n(n+1)}{2}\right)^{5}-b_{2}\cdot\left(\frac{n(n+1)}{2}\right)^{4}+b_{3}\cdot\left(\frac{n(n+1)}{2}\right)^{3}-b_{4}\cdot\left(\frac{n(n+1)}{2}\right)^{2}+b_{5}\cdot O_{5}.\label{eq:73}
\end{equation}
 The system of equations is 
\begin{align*}
O_{13}: &  & -b_{2}\cdot\left(\frac{n(n+1)}{2}\right)^{4} & +b_{3}\cdot\left(\frac{n(n+1)}{2}\right)^{3} & -b_{4}\cdot\left(\frac{n(n+1)}{2}\right)^{2} & +b_{5}\cdot O_{5}\\
O_{11}: &  & 35\cdot\frac{32}{6}\cdot\left(\frac{n(n+1)}{2}\right)^{4} & -35\cdot\frac{64}{6}\cdot\left(\frac{n(n+1)}{2}\right)^{3} & +35\cdot\frac{68}{6}\cdot\left(\frac{n(n+1)}{2}\right)^{2} & -35\cdot\frac{30}{6}\cdot O_{5}\\
O_{9}: &  &  & 21\cdot\frac{16}{5}\cdot\left(\frac{n(n+1)}{2}\right)^{3} & -21\cdot\frac{20}{5}\cdot\left(\frac{n(n+1)}{2}\right)^{2} & +21\cdot\frac{9}{5}\cdot O_{5}\\
O_{7}: &  &  &  & 1\cdot\frac{8}{4}\cdot\left(\frac{n(n+1)}{2}\right)^{2} & -1\cdot\frac{4}{4}\cdot O_{5},
\end{align*}
 where 
\begin{eqnarray*}
 & -b_{2}+35\cdot\frac{32}{6}=0\\
 & b_{3}-35\cdot\frac{64}{6}+21\cdot\frac{16}{5}=0\\
 & -b_{4}+35\cdot\frac{68}{6}-21\cdot\frac{20}{5}+1\cdot\frac{8}{4}=0\\
 & b_{5}-35\cdot\frac{30}{6}+21\cdot\frac{9}{5}-1\cdot\frac{4}{4}=0.
\end{eqnarray*}
 We can solve it immediately: 
\begin{align*}
b_{2} & =\frac{560}{3}\\
b_{3} & =\frac{4592}{15}\\
b_{4} & =\frac{944}{3}\\
b_{5} & =\frac{691}{5}.
\end{align*}
 Therefore we can rewrite expression \ref{eq:73} as 
\begin{equation}
7\cdot O_{13}=64\cdot\left(\frac{n(n+1)}{2}\right)^{5}-\frac{560}{3}\cdot\left(\frac{n(n+1)}{2}\right)^{4}+\frac{4592}{15}\cdot\left(\frac{n(n+1)}{2}\right)^{3}-\frac{944}{3}\cdot\left(\frac{n(n+1)}{2}\right)^{2}+\frac{691}{5}\cdot O_{5}.
\end{equation}
 If we insert it into 
\[
\sum k^{13}=O_{13}\cdot\sum k^{3}
\]
 then we can test the result: 
\[
\sum_{k=1}^{9}k^{13}=3,202,860,761,145
\]
 and 
\[
\frac{1}{7}\cdot\left(64\cdot45^{5}-\frac{560}{3}\cdot45^{4}+\frac{4592}{15}\cdot45^{3}-\frac{944}{3}\cdot45^{2}+\frac{691}{5}\cdot\frac{1}{3}\cdot\left(4\cdot45-1\right)\right)\cdot\sum_{k=1}^{9}k^{3}
\]
 
\[
=\frac{1}{7}\cdot\frac{166074261689}{15}\cdot2025=3,202,860,761,145.
\]

\section{A Long List 2}

One last time, let us place our results for $E_{2m}$ and $O_{2m+1}$
into a list: 
\begin{align*}
\\
E_{2} & =1\\
5\cdot E_{4} & =6\cdot\left(\frac{n(n+1)}{2}\right)-1\\
7\cdot E_{6} & =12\cdot\left(\frac{n(n+1)}{2}\right)^{2}-5\cdot E_{4}\\
9\cdot E_{8} & =24\cdot\left(\frac{n(n+1)}{2}\right)^{3}-24\cdot\left(\frac{n(n+1)}{2}\right)^{2}+9\cdot E_{4}\\
11\cdot E_{10} & =48\cdot\left(\frac{n(n+1)}{2}\right)^{4}-80\cdot\left(\frac{n(n+1)}{2}\right)^{3}+68\cdot\left(\frac{n(n+1)}{2}\right)^{2}-25\cdot E_{4}\\
13\cdot E_{12} & =96\cdot\left(\frac{n(n+1)}{2}\right)^{5}-240\cdot\left(\frac{n(n+1)}{2}\right)^{4}+328\cdot\left(\frac{n(n+1)}{2}\right)^{3}-\frac{1888}{7}\cdot\left(\frac{n(n+1)}{2}\right)^{2}+\frac{691}{7}\cdot E_{4}\\
\\
O_{3} & =1\\
3\cdot O_{5} & =4\cdot\left(\frac{n(n+1)}{2}\right)-1\\
4\cdot O_{7} & =8\cdot\left(\frac{n(n+1)}{2}\right)^{2}-4\cdot O_{5}\\
5\cdot O_{9} & =16\cdot\left(\frac{n(n+1)}{2}\right)^{3}-20\cdot\left(\frac{n(n+1)}{2}\right)^{2}+9\cdot O_{5}\\
6\cdot O_{11} & =32\cdot\left(\frac{n(n+1)}{2}\right)^{4}-64\cdot\left(\frac{n(n+1)}{2}\right)^{3}+68\cdot\left(\frac{n(n+1)}{2}\right)^{2}-30\cdot O_{5}\\
7\cdot O_{13} & =64\cdot\left(\frac{n(n+1)}{2}\right)^{5}-\frac{560}{3}\cdot\left(\frac{n(n+1)}{2}\right)^{4}+\frac{4592}{15}\cdot\left(\frac{n(n+1)}{2}\right)^{3}-\frac{944}{3}\cdot\left(\frac{n(n+1)}{2}\right)^{2}+\frac{691}{5}\cdot O_{5}.\\
\end{align*}

As a reminder, the expressions for $E_{10},E_{12}$ and $O_{11},O_{13}$
are conjectures only. Do we have any observations?
\begin{enumerate}
\item we have come a long way. The list is a large improvement over the
one in Section \ref{sec:A-Long-List}. The leading fractions, the
coefficients for $\frac{n(n+1)}{2}$, the alternating signs: everything
is in unison. We truly believe we have found the correct expressions.
Of course we doubt we could have guessed coefficients like 
\[
\frac{1888}{7},\frac{691}{7};\,\frac{560}{3},\frac{4592}{15},\frac{944}{3},\frac{691}{5},
\]
 and we find their appearance even a bit strange, but now we can explain
where they come from.
\end{enumerate}
\pagebreak{}

\section{Summary of Part 4}

We began Part 4 with our most promising results from Part 3, the lists
\begin{align*}
1 & =1\\
2\cdot3 & =1+5\\
2^{2}\cdot3 & =0+5+7\\
2^{3}\cdot3 & =0+1+14+9\\
2^{4}\cdot3 & =0+0+7+30+11
\end{align*}
 and 
\begin{align*}
1 & =1\\
2^{2} & =1+3\\
2^{3} & =0+4+4\\
2^{4} & =0+1+10+5\\
2^{5} & =0+0+6+20+6,
\end{align*}
 which contained relationships for terms in the expressions for $E_{2m}$
and $O_{2m+1}$. The outward appearances of the lists suggested we
turn to Pascal's Triangle. In little short of a miracle we noticed
\begin{align*}
1 & =1\\
2\cdot3 & =\left[\binom{2}{-1}+\binom{3}{0}\right]+\left[\binom{2}{1}+\binom{3}{2}\right]\\
2^{2}\cdot3 & =\left[\binom{3}{0}+\binom{4}{1}\right]+\left[\binom{3}{2}+\binom{4}{3}\right]\\
2^{3}\cdot3 & =\left[\binom{4}{-1}+\binom{5}{0}\right]+\left[\binom{4}{1}+\binom{5}{2}\right]+\left[\binom{4}{3}+\binom{5}{4}\right]\\
2^{4}\cdot3 & =\left[\binom{5}{0}+\binom{6}{1}\right]+\left[\binom{5}{2}+\binom{6}{3}\right]+\left[\binom{5}{4}+\binom{6}{5}\right]
\end{align*}
 and 
\begin{align*}
1 & =1\\
2^{2} & =\binom{3}{0}+\binom{3}{2}\\
2^{3} & =\binom{4}{1}+\binom{4}{3}\\
2^{4} & =\binom{5}{0}+\binom{5}{2}+\binom{5}{4}\\
2^{5} & =\binom{6}{1}+\binom{6}{3}+\binom{6}{5}.
\end{align*}
 All of our previous difficulties seemed to disappear. After guessing
expressions for the general terms, 
\begin{align*}
2^{m-1}\cdot3 & =\left[\binom{m}{-1}+\binom{m+1}{0}\right]+\left[\binom{m}{1}+\binom{m+1}{2}\right]+\left[\binom{m}{3}+\binom{m+1}{4}\right]\\
 & +\cdots+\left[\binom{m}{m-1}+\binom{m+1}{m}\right]\\
2^{m}\cdot3 & =\left[\binom{m+1}{0}+\binom{m+2}{1}\right]+\left[\binom{m+1}{2}+\binom{m+2}{3}\right]+\left[\binom{m+1}{4}+\binom{m+2}{5}\right]\\
 & +\cdots+\left[\binom{m+1}{m}+\binom{m+2}{m+1}\right]
\end{align*}
 and 
\begin{align*}
2^{m} & =\binom{m+1}{0}+\binom{m+1}{2}+\binom{m+1}{4}+\cdots+\binom{m+1}{m}\\
2^{m+1} & =\binom{m+2}{1}+\binom{m+2}{3}+\binom{m+2}{5}+\cdots+\binom{m+2}{m+1},
\end{align*}
 and recasting some previous relationships in terms of entries in
Pascal's Triangle, we used the new results to calculate the next cases
of $E_{12}$ and $O_{13}$. With less effort than before we found
\[
13\cdot E_{12}=96\cdot\left(\frac{n(n+1)}{2}\right)^{5}-240\cdot\left(\frac{n(n+1)}{2}\right)^{4}+328\cdot\left(\frac{n(n+1)}{2}\right)^{3}-\frac{1888}{7}\cdot\left(\frac{n(n+1)}{2}\right)^{2}+\frac{691}{7}\cdot E_{4}
\]
\[
7\cdot O_{13}=64\cdot\left(\frac{n(n+1)}{2}\right)^{5}-\frac{560}{3}\cdot\left(\frac{n(n+1)}{2}\right)^{4}+\frac{4592}{15}\cdot\left(\frac{n(n+1)}{2}\right)^{3}-\frac{944}{3}\cdot\left(\frac{n(n+1)}{2}\right)^{2}+\frac{691}{5}\cdot O_{5},
\]
 which we were able to verify in the special cases of 
\[
\sum_{k=1}^{9}k^{12}=E_{12}(n)\cdot\sum_{k=1}^{9}k^{2}
\]
\[
\sum_{k=1}^{9}k^{13}=O_{13}(n)\cdot\sum_{k=1}^{9}k^{3}.
\]
 Finally, we placed all of our expressions for $E_{2m}$ and $O_{2m+1}$
into a long list and made some quick observations.

\pagebreak{}

\section{Final Remarks\label{sec:Final-Remarks}}

At this point we believe we have answered the question we posed at
the start,
\begin{quotation}
``What is a general expression for 
\[
\sum_{k=1}^{n}k^{m}=1^{m}+2^{m}+3^{m}+\cdots+n^{m},
\]
where $m$ is a positive integer?''
\end{quotation}
The answer is,
\begin{quotation}
``For even powers the expression is 
\begin{eqnarray*}
\left(2m+1\right)\cdot\frac{\sum_{k=1}^{n}k^{2m}}{\sum_{k=1}^{n}k^{2}} & = & a_{1}\cdot\left(\frac{n(n+1)}{2}\right)^{m-1}-a_{2}\cdot\left(\frac{n(n+1)}{2}\right)^{m-2}+a_{3}\cdot\left(\frac{n(n+1)}{2}\right)^{m-3}\\
 &  & \mp\cdots\mp a_{m-2}\cdot\left(\frac{n(n+1)}{2}\right)^{2}\mp a_{m-1}\cdot\frac{1}{5}\cdot\left(6\cdot\frac{n(n+1)}{2}-1\right),
\end{eqnarray*}
 and for odd powers, 
\begin{eqnarray*}
\left(m+1\right)\cdot\frac{\sum_{k=1}^{n}k^{2m+1}}{\sum_{k=1}^{n}k^{3}} & = & b_{1}\cdot\left(\frac{n(n+1)}{2}\right)^{m-1}-b_{2}\cdot\left(\frac{n(n+1)}{2}\right)^{m-2}+b_{3}\cdot\left(\frac{n(n+1)}{2}\right)^{m-3}\\
 &  & \mp\cdots\mp b_{m-2}\cdot\left(\frac{n(n+1)}{2}\right)^{2}\mp b_{m-1}\cdot\frac{1}{3}\cdot\left(4\cdot\frac{n(n+1)}{2}-1\right),
\end{eqnarray*}
 where $a_{i}$ and $b_{j}$ are rational numbers.''
\end{quotation}
Whether or not we were aware of it previously, what our work has produced
is an algorithm for finding expressions for $\sum k^{m}$. Even though
it still is a conjecture, it is far preferable to the established
result of Section \ref{sub:3 The First Catch}, 
\[
\sum_{k=1}^{n}k^{m+1}+\sum_{k=1}^{n}\sum_{l=1}^{k}l^{m}=\left(n+1\right)\cdot\sum_{k=1}^{n}k^{m},
\]
 which involves wading through lengthy calculations. With the old
result we gave up on calculations after $\sum k^{9}$. A case as far
out as $\sum k^{16}$ almost would be unthinkable. It would require
taming an expression like 
\[
\sum\sum l^{15}=\frac{c_{1}\cdot\sum k+c_{2}\cdot\sum k^{2}+c_{3}\cdot\sum k^{3}+\cdots+c_{16}\cdot\sum k^{16}}{C},
\]
 which might contain as many as nine terms. For comparison, with the
new algorithm we can do the following.

\pagebreak{}

We start with 
\[
\sum k^{16}=E_{16}\cdot\sum k^{2}.
\]
 The exact form of $E_{16}$ satisfies 
\begin{equation}
e_{2}E_{2}+e_{4}E_{4}+e_{6}E_{6}+e_{8}E_{8}+e_{10}E_{10}+e_{12}E_{12}+e_{14}E_{14}+e_{16}E_{16}=384\cdot\left(\frac{n(n+1)}{2}\right)^{7}.\label{eq:75}
\end{equation}
 In order to find the values for $e_{2},e_{4},e_{6},\ldots,e_{16}$
we look at the previous list of 
\begin{align*}
1 & =1\\
6 & =1+5\\
12 & =0+5+7\\
24 & =0+1+14+9\\
48 & =0+0+7+30+11\\
96 & =0+0+1+27+55+13\\
192 & =0+0+0+9+77+91+15.
\end{align*}
 It does not contain the values for $E_{16}$. Not to worry. From
the relationships in Pascal's Triangle we notice 
\begin{align*}
2^{5}\cdot3 & =\left[\binom{6}{-1}+\binom{7}{0}\right]+\left[\binom{6}{1}+\binom{7}{2}\right]+\left[\binom{6}{3}+\binom{7}{4}\right]+\left[\binom{6}{5}+\binom{7}{6}\right]\\
2^{6}\cdot3 & =\left[\binom{7}{0}+\binom{8}{1}\right]+\left[\binom{7}{2}+\binom{8}{3}\right]+\left[\binom{7}{4}+\binom{8}{5}\right]+\left[\binom{7}{6}+\binom{8}{7}\right],
\end{align*}
 which suggest the next entry will be 
\begin{eqnarray*}
2^{7}\cdot3 & = & \left[\binom{8}{-1}+\binom{9}{0}\right]+\left[\binom{8}{1}+\binom{9}{2}\right]+\left[\binom{8}{3}+\binom{9}{4}\right]+\left[\binom{8}{5}+\binom{9}{6}\right]+\left[\binom{8}{7}+\binom{9}{8}\right]\\
 &  & =1+44+182+140+17.
\end{eqnarray*}
 We write it as 
\[
384=0+0+0+1+44+182+140+17,
\]
 which tells us $e_{2}=e_{4}=e_{6}=0$. We substitute these values
into expression \ref{eq:75}, 
\[
1\cdot E_{8}+44\cdot E_{10}+182\cdot E_{12}+140\cdot E_{14}+17\cdot E_{16}=384\cdot\left(\frac{n(n+1)}{2}\right)^{7},
\]
 and then simplify 
\[
17\cdot E_{16}=384\cdot\left(\frac{n(n+1)}{2}\right)^{7}-\left(E_{8}+44\cdot E_{10}+182\cdot E_{12}+140\cdot E_{14}\right),
\]
 which, for comparison, contains only four terms. The final expression
will be 
\begin{eqnarray*}
17\cdot E_{16} & = & 384\cdot\left(\frac{n(n+1)}{2}\right)^{7}-a_{2}\cdot\left(\frac{n(n+1)}{2}\right)^{6}+a_{3}\cdot\left(\frac{n(n+1)}{2}\right)^{5}\\
 &  & \mp\cdots-a_{6}\cdot\left(\frac{n(n+1)}{2}\right)^{2}+a_{7}\cdot\frac{1}{5}\cdot\left(6\cdot\frac{n(n+1)}{2}-1\right)
\end{eqnarray*}
 for some rational numbers $a_{2},a_{3},\ldots,a_{6},a_{7}$. Last,
we reinsert it into 
\[
\sum k^{16}=E_{16}\cdot\sum k^{2}.
\]

Unfortunately, the algorithm still is recursive. In order to find
the expression for $E_{16}$ we need to know already the expressions
for $E_{2},E_{4},E_{6},\ldots,E_{14}$. Can we improve upon this?
For the previous expressions for $E_{2m}$ and $O_{2m+1}$ it is tempting
to replace the coefficients with the entries from Pascal's Triangle
and then to try to rewrite them into explicit expressions. However,
to proceed that way seems to place us back into lengthy calculations
reminiscent of the old result. Therefore, for the question at the
start of the paper we have an answer, but we cannot give the exact
values for the coefficients, $a_{i}$ and $b_{j}$. We have to find
them.

\bigskip{}

\rule[0.5ex]{0.5\columnwidth}{1pt}

\bigskip{}

\begin{description}
\item [{Historical~Note}] After the author had completed the paper he
was informed he had rediscovered a result credited to Johann Faulhaber,
1580-1635. For a discussion of Faulhaber's original work and some
developments over the ensuing centuries, the reader can consult the
paper ``Johann Faulhaber and Sums of Powers'' by Donald E. Knuth.
Perhaps at a later time the author will address the matter himself.
\end{description}
\pagebreak{}

\section*{Acknowledgements}

The author expresses his gratitude to George Polya, especially for
the works ``How to Solve It'' and ``Mathematics and Plausible Reasoning.''
The author expresses his gratitude also to the many people who contribute
to free and open source software, especially Slackware Linux, \LyX{},
and Maxima.

On a lighter note, the author thanks Mark Twain for writing ``The
Adventures of Huckleberry Finn,'' which he had been reading in his
time away from working on this paper and which had provided the inspiration
for some of the section titles.

Last, the author is indebted to someone who wishes to remain anonymous,
who was kind enough to read a rough version of the paper and suggested
ways to improve it.

\pagebreak{}

\appendix

\section{Appendix\label{sec:Appendix}}

\subsection{$\sum_{k=1}^{n}k^{7}$\label{sub:seven}}

What is 
\[
\sum_{k=1}^{n}k^{7}=1^{7}+2^{7}+3^{7}+\cdots+n^{7}\,?
\]
 We know that 
\begin{equation}
\sum_{k=1}^{n}k^{7}=\left(n+1\right)\cdot\sum_{k=1}^{n}k^{6}-\sum_{k=1}^{n}\sum_{l=1}^{k}l^{6}.\label{eq:76}
\end{equation}
 We know also that 
\[
\sum_{k=1}^{n}k^{6}=\frac{n-7n^{3}+21n^{5}+21n^{6}+6n^{7}}{42},
\]
 which implies 
\[
\sum_{k=1}^{n}\sum_{l=1}^{k}l^{6}=\frac{\sum_{k=1}^{n}k-7\cdot\sum_{k=1}^{n}k^{3}+21\cdot\sum_{k=1}^{n}k^{5}+21\cdot\sum_{k=1}^{n}k^{6}+6\cdot\sum_{k=1}^{n}k^{7}}{42}.
\]
 Therefore we may rewrite expression \ref{eq:76} as 
\begin{eqnarray*}
\sum_{k=1}^{n}k^{7} & = & \left(n+1\right)\cdot\sum_{k=1}^{n}k^{6}-\frac{\sum_{k=1}^{n}k-7\cdot\sum_{k=1}^{n}k^{3}+21\cdot\sum_{k=1}^{n}k^{5}+21\cdot\sum_{k=1}^{n}k^{6}+6\cdot\sum_{k=1}^{n}k^{7}}{42}\\
 &  & =\left(n+1\right)\cdot\sum_{k=1}^{n}k^{6}-\frac{1}{42}\cdot\sum_{k=1}^{n}k+\frac{1}{6}\cdot\sum_{k=1}^{n}k^{3}-\frac{1}{2}\cdot\sum_{k=1}^{n}k^{5}-\frac{1}{2}\cdot\sum_{k=1}^{n}k^{6}-\frac{1}{7}\cdot\sum_{k=1}^{n}k^{7}
\end{eqnarray*}
\[
\frac{8}{7}\cdot\sum_{k=1}^{n}k^{7}=\frac{2n+1}{2}\cdot\sum_{k=1}^{n}k^{6}-\frac{1}{42}\cdot\sum_{k=1}^{n}k+\frac{1}{6}\cdot\sum_{k=1}^{n}k^{3}-\frac{1}{2}\cdot\sum_{k=1}^{n}k^{5}.
\]

\bigskip{}
For $\sum_{k=1}^{n}k^{5}$ we substitute 
\[
\frac{2n\left(n+1\right)-1}{3}\cdot\sum_{k=1}^{n}k^{3}
\]
 and collect the terms which involve $\sum_{k=1}^{n}k^{3}$: 
\[
\frac{8}{7}\cdot\sum_{k=1}^{n}k^{7}=\frac{2n+1}{2}\cdot\sum_{k=1}^{n}k^{6}-\frac{1}{42}\cdot\sum_{k=1}^{n}k+\left(\frac{1}{6}-\frac{1}{2}\cdot\frac{2n\left(n+1\right)-1}{3}\right)\cdot\sum_{k=1}^{n}k^{3}
\]
\[
\frac{8}{7}\cdot\sum_{k=1}^{n}k^{7}=\frac{2n+1}{2}\cdot\sum_{k=1}^{n}k^{6}-\frac{1}{42}\cdot\sum_{k=1}^{n}k+\left(\frac{1-n\left(n+1\right)}{3}\right)\cdot\sum_{k=1}^{n}k^{3}.
\]
 Since 
\[
\sum_{k=1}^{n}k^{3}=\left(\frac{n\left(n+1\right)}{2}\right)^{2}=\sum_{k=1}^{n}k\cdot\sum_{k=1}^{n}k,
\]
 we collect terms of the form $\sum_{k=1}^{n}k$: 
\[
\frac{8}{7}\cdot\sum_{k=1}^{n}k^{7}=\frac{2n+1}{2}\cdot\sum_{k=1}^{n}k^{6}-\frac{1}{42}\cdot\sum_{k=1}^{n}k+\left(\frac{1-n\left(n+1\right)}{3}\right)\cdot\sum_{k=1}^{n}k\cdot\sum_{k=1}^{n}k
\]
\[
\frac{8}{7}\cdot\sum_{k=1}^{n}k^{7}=\frac{2n+1}{2}\cdot\sum_{k=1}^{n}k^{6}-\left(\frac{1}{42}-\left(\frac{1-n\left(n+1\right)}{3}\right)\cdot\frac{n\left(n+1\right)}{2}\right)\cdot\sum_{k=1}^{n}k
\]
\[
\frac{8}{7}\cdot\sum_{k=1}^{n}k^{7}=\frac{2n+1}{2}\cdot\sum_{k=1}^{n}k^{6}-\left(\frac{1}{42}-\frac{n\left(n+1\right)-\left(n\left(n+1\right)\right)^{2}}{6}\right)\cdot\sum_{k=1}^{n}k
\]
\[
\frac{8}{7}\cdot\sum_{k=1}^{n}k^{7}=\frac{2n+1}{2}\cdot\sum_{k=1}^{n}k^{6}-\left(\frac{7\left(n\left(n+1\right)\right)^{2}-7n\left(n+1\right)+1}{42}\right)\cdot\sum_{k=1}^{n}k,
\]
 which allows us to write 
\begin{equation}
\sum_{k=1}^{n}k^{7}=\frac{7}{8}\cdot\frac{2n+1}{2}\cdot\sum_{k=1}^{n}k^{6}-\left(\frac{7\left(n\left(n+1\right)\right)^{2}-7n\left(n+1\right)+1}{48}\right)\cdot\sum_{k=1}^{n}k.\label{eq:77}
\end{equation}

For $\sum_{k=1}^{n}k^{6}$ we substitute 
\[
\sum_{k=1}^{n}k^{6}=\left(\frac{3\left(n\left(n+1\right)\right)^{2}-3n\left(n+1\right)+1}{7}\right)\cdot\sum_{k=1}^{n}k^{2}
\]
 and rewrite the left terms of expression \ref{eq:77} as 
\[
\frac{7}{8}\cdot\frac{2n+1}{2}\cdot\left(\frac{3\left(n\left(n+1\right)\right)^{2}-3n\left(n+1\right)+1}{7}\right)\cdot\sum_{k=1}^{n}k^{2}
\]
\[
=\frac{7}{8}\cdot\frac{2n+1}{2}\cdot\left(\frac{3\left(n\left(n+1\right)\right)^{2}-3n\left(n+1\right)+1}{7}\right)\cdot\frac{2n+1}{3}\cdot\frac{n\left(n+1\right)}{2}
\]
\[
=\frac{7}{8}\cdot\frac{2n+1}{2}\cdot\frac{2n+1}{3}\cdot\left(\frac{3\left(n\left(n+1\right)\right)^{2}-3n\left(n+1\right)+1}{7}\right)\cdot\frac{n\left(n+1\right)}{2}
\]
\[
=\frac{7\left(2n+1\right)^{2}}{48}\cdot\left(\frac{3\left(n\left(n+1\right)\right)^{2}-3n\left(n+1\right)+1}{7}\right)\cdot\sum_{k=1}^{n}k
\]
\[
=\left(2n+1\right)^{2}\cdot\left(\frac{3\left(n\left(n+1\right)\right)^{2}-3n\left(n+1\right)+1}{48}\right)\cdot\sum_{k=1}^{n}k.
\]
 We substitute this back into expression \ref{eq:77} to get 
\[
\sum_{k=1}^{n}k^{7}=\left(2n+1\right)^{2}\cdot\left(\frac{3\left(n\left(n+1\right)\right)^{2}-3n\left(n+1\right)+1}{48}\right)\cdot\sum_{k=1}^{n}k
\]
\[
-\left(\frac{7\left(n\left(n+1\right)\right)^{2}-7n\left(n+1\right)+1}{48}\right)\cdot\sum_{k=1}^{n}k
\]
\[
=\left(\left(2n+1\right)^{2}\cdot\frac{3\left(n\left(n+1\right)\right)^{2}-3n\left(n+1\right)+1}{48}-\frac{7\left(n\left(n+1\right)\right)^{2}-7n\left(n+1\right)+1}{48}\right)\cdot\sum_{k=1}^{n}k
\]
\begin{equation}
=\frac{1}{8}\cdot\left(\left(2n+1\right)^{2}\cdot\frac{3\left(n\left(n+1\right)\right)^{2}-3n\left(n+1\right)+1}{6}-\frac{7\left(n\left(n+1\right)\right)^{2}-7n\left(n+1\right)+1}{6}\right)\cdot\sum_{k=1}^{n}k.\label{eq:78}
\end{equation}
 Consider expression \ref{eq:78} a first derivation for $\sum_{k=1}^{n}k^{7}$.

\bigskip{}
Let us continue. We rewrite the expression inside the parentheses
as 
\[
\left(2n+1\right)^{2}\cdot2\cdot\left(\frac{n\left(n+1\right)}{2}\right)^{2}-\left(2n+1\right)^{2}\cdot\frac{n\left(n+1\right)}{2}+\left(2n+1\right)^{2}\cdot\frac{1}{6}
\]
\[
-\frac{14}{3}\cdot\left(\frac{n\left(n+1\right)}{2}\right)^{2}+\frac{7}{3}\cdot\frac{n\left(n+1\right)}{2}-\frac{1}{6}
\]
\[
=\left(2\cdot\left(2n+1\right)^{2}-\frac{14}{3}\right)\cdot\left(\frac{n\left(n+1\right)}{2}\right)^{2}+\left(-\left(2n+1\right)^{2}+\frac{7}{3}\right)\cdot\frac{n\left(n+1\right)}{2}+\left(\frac{\left(2n+1\right)^{2}}{6}-\frac{1}{6}\right)
\]
\[
=\frac{24n^{2}+24n-8}{3}\cdot\left(\frac{n\left(n+1\right)}{2}\right)^{2}-\frac{12n^{2}+12n-4}{3}\cdot\frac{n\left(n+1\right)}{2}+\frac{4}{3}\cdot\frac{n\left(n+1\right)}{2}.
\]
 If we multiply the expression by 
\[
\sum_{k=1}^{n}k=\frac{n\left(n+1\right)}{2}
\]
 then we get 
\[
\frac{24n^{2}+24n-8}{3}\cdot\frac{n\left(n+1\right)}{2}\cdot\left(\frac{n\left(n+1\right)}{2}\right)^{2}-\frac{12n^{2}+12n-4}{3}\cdot\frac{n\left(n+1\right)}{2}\cdot\frac{n\left(n+1\right)}{2}
\]
\[
+\frac{4}{3}\cdot\frac{n\left(n+1\right)}{2}\cdot\frac{n\left(n+1\right)}{2}
\]
\[
=\frac{24n^{2}+24n-8}{3}\cdot\frac{n\left(n+1\right)}{2}\cdot\sum_{k=1}^{n}k^{3}-\frac{12n^{2}+12n-4}{3}\cdot\sum_{k=1}^{n}k^{3}+\frac{4}{3}\cdot\sum_{k=1}^{n}k^{3}
\]
\[
=\frac{8\left(3(n(n+1))^{2}-n(n+1)\right)}{6}\cdot\sum_{k=1}^{n}k^{3}-\frac{4\left(3n(n+1)-1\right)}{3}\cdot\sum_{k=1}^{n}k^{3}+\frac{4}{3}\cdot\sum_{k=1}^{n}k^{3}.
\]
 Remembering the form of expression \ref{eq:78}, if we multiply by
$\frac{1}{8}$ then we get 
\[
\sum_{k=1}^{n}k^{7}=\frac{1}{8}\cdot\left(\frac{8\left(3(n(n+1))^{2}-n(n+1)\right)}{6}\cdot\sum_{k=1}^{n}k^{3}-\frac{4\left(3n(n+1)-1\right)}{3}\cdot\sum_{k=1}^{n}k^{3}+\frac{4}{3}\cdot\sum_{k=1}^{n}k^{3}\right)
\]
\[
=\frac{1}{8}\cdot\left(\frac{8\left(3(n(n+1))^{2}-n(n+1)\right)}{6}-\frac{4\left(3n(n+1)-1\right)}{3}+\frac{4}{3}\right)\cdot\sum_{k=1}^{n}k^{3}
\]
\[
=\frac{1}{8}\cdot\left(\frac{12(n(n+1))^{2}-16n(n+1)+8}{3}\right)\cdot\sum_{k=1}^{n}k^{3}
\]
\begin{equation}
=\left(\frac{3(n(n+1))^{2}-4n(n+1)+2}{6}\right)\cdot\sum_{k=1}^{n}k^{3}.\label{eq:79}
\end{equation}
 Another way to write it is 
\[
\sum_{k=1}^{n}k^{7}=\frac{1}{8}\cdot\left(\frac{12(n(n+1))^{2}-16n(n+1)+8}{3}\right)\cdot\sum_{k=1}^{n}k^{3}
\]
\[
=\frac{1}{8}\cdot\left(\frac{48}{3}\cdot\left(\frac{n(n+1)}{2}\right)^{2}-\frac{32}{3}\cdot\frac{n(n+1)}{2}+\frac{8}{3}\right)\cdot\sum_{k=1}^{n}k^{3}
\]
\begin{equation}
=\frac{1}{3}\cdot\left(6\cdot\left(\frac{n(n+1)}{2}\right)^{2}-4\cdot\frac{n(n+1)}{2}+1\right)\cdot\sum_{k=1}^{n}k^{3}.\label{eq:80}
\end{equation}
 Both expressions simplify to 
\[
\sum_{k=1}^{n}k^{7}=\frac{2n^{2}-7n^{4}+14n^{6}+12n^{7}+3n^{8}}{24}.
\]

\pagebreak{}

\subsection{Does a prime $p=2m+1$ divide $\sum_{k=1}^{m}k^{2}$?\label{sub:Does--divide}}

In the calculations for $E_{10}$, under the assumption that both
$e_{8}$ and $a_{2}$ are integers, the equation 
\[
\frac{8}{3}\cdot e_{8}=a_{2}
\]
 tells us that $e_{8}=3c$ for some integer $c$. If we were to carry
out the calculations for $E_{12}$, the analogous conjectures would
be 
\[
e_{4}E_{4}+e_{6}E_{6}+e_{8}E_{8}+e_{10}E_{10}+13E_{12}=96\cdot\left(\frac{n(n+1)}{2}\right)^{5}
\]
 and 
\[
13\cdot E_{12}=96\cdot\left(\frac{n(n+1)}{2}\right)^{5}-a_{2}^{*}\cdot\left(\frac{n(n+1)}{2}\right)^{4}+a_{3}\cdot\left(\frac{n(n+1)}{2}\right)^{3}-a_{4}\cdot\left(\frac{n(n+1)}{2}\right)^{2}+a_{5}\cdot E_{4}
\]
 and the list of coefficients would be 
\begin{align*}
1 & =1\\
6 & =1+5\\
12 & =0+5+7\\
24 & =0+1+14+9\\
48 & =0+0+7+30+11\\
96 & =0+e_{4}+e_{6}+e_{8}+e_{10}+13.
\end{align*}
 Looking at the diagonal containing $e_{10}$, we would suspect $e_{10}=55=\sum_{k=1}^{5}k^{2}$.

Further, under an analogous system of equations, under the assumption
that both $e_{10}$ and $a_{2}^{*}$ were integers, the equation 
\[
\frac{48}{11}\cdot e_{10}=a_{2}^{*}
\]
 would tell us $e_{10}=11d$ for some integer $d.$ Therefore we would
have $11d=55=\sum_{k=1}^{5}k^{2}$, which means 11 divides $55=\sum_{k=1}^{5}k^{2}$.
We write $11\mid55=\sum_{k=1}^{5}k^{2}.$ That happens to be true:
11$\cdot$5=55. If we look again at the list of coefficients then
we notice 
\[
3\mid30=\sum_{k=1}^{4}k^{2},\,7\mid14=\sum_{k=1}^{3}k^{2},\,5\mid5=\sum_{k=1}^{2}k^{2}.
\]
 If we ignore the case of $3\mid30=\sum_{k=1}^{4}k^{2}$ for the moment,
we may make the following observation:
\begin{description}
\item [{Exercise}] if $p=2m+1$ is prime then $p$ divides $\sum_{k=1}^{m}k^{2}$.
Prove or disprove.
\end{description}
\pagebreak{}

\subsection{Inductive Reasoning Sometimes Goes Wrong\label{sub:Inductive-Reasoning-Sometimes}}

In Section \ref{sub:3 The First Catch} we proved formally that 
\[
\sum_{k=1}^{n}k^{m+1}+\sum_{k=1}^{n}\sum_{l=1}^{k}l^{m}=\left(n+1\right)\cdot\sum_{k=1}^{n}k^{m}
\]
 for positive integers $m$. Over the course of Parts 1 and 2 we used
it to derive 
\[
\sum k^{3}=O_{3}\cdot\sum k^{3}
\]
\[
\sum k^{5}=O_{5}\cdot\sum k^{3}
\]
\[
\sum k^{7}=O_{7}\cdot\sum k^{3}
\]
\[
\sum k^{9}=O_{9}\cdot\sum k^{3},
\]
 where 
\[
O_{3}=1
\]
\[
O_{5}=\frac{1}{3}\cdot\left(4\cdot\frac{n\left(n+1\right)}{2}-1\right)
\]

\[
O_{7}=\frac{1}{4}\cdot\left(8\cdot\left(\frac{n(n+1)}{2}\right)^{2}-4\cdot O_{5}\right)
\]

\begin{equation}
O_{9}=\frac{1}{5}\cdot\left(16\cdot\left(\frac{n(n+1)}{2}\right)^{3}-20\cdot\left(\frac{n(n+1)}{2}\right)^{2}+9\cdot O_{5}\right).\label{eq:81}
\end{equation}
 The expressions for $\sum k^{3},\sum k^{5},\sum k^{7},$ and $\sum k^{9}$
were true. In Part 3 we used inductive reasoning to arrive at 
\begin{equation}
O_{11}=\frac{1}{6}\cdot\left(32\cdot\left(\frac{n(n+1)}{2}\right)^{4}-64\cdot\left(\frac{n(n+1)}{2}\right)^{3}+68\cdot\left(\frac{n(n+1)}{2}\right)^{2}-30\cdot O_{5}\right),\label{eq:82}
\end{equation}
 which we believed satisfied 
\[
\sum k^{11}=O_{11}\cdot\sum k^{3}.
\]
 As a check on our inductive reasoning we verified it in the special
case of 
\[
\sum_{k=1}^{9}k^{11}=42,364,319,625.
\]
 Given how complicated the sums were, we believed that was sufficient
evidence to suggest the result was correct. What if instead we had
verified it in the following way?

In the expressions in (\ref{eq:81}) the coefficients satisfy 
\begin{align*}
3 & =4-1\\
4 & =8-4\\
5 & =16-20+9.
\end{align*}
 In the new result for $O_{11}$ they satisfy 
\begin{align*}
6 & =32-64+68-30.
\end{align*}
 Impressive, right? It also bears an analogy with our previous conjectures.
Is it enough to justify we found the correct expression? Not quite.

\bigskip{}
In the derivation for $O_{11}$ we had the conjectures 
\[
o_{5}O_{5}+o_{7}O_{7}+o_{9}O_{9}+o_{11}O_{11}=32\cdot\left(\frac{n(n+1)}{2}\right)^{4}
\]
 and 
\begin{equation}
o_{5}+o_{7}+o_{9}+o_{11}=32\label{eq:83}
\end{equation}
 and the system of equations 
\begin{eqnarray*}
 & -b_{2}+o_{9}\cdot\frac{32}{10}=0\\
 & b_{3}-o_{9}\cdot\frac{40}{10}+o_{7}\cdot\frac{6}{3}=0\\
 & -b_{4}+o_{9}\cdot\frac{18}{10}-o_{7}\cdot\frac{3}{3}+o_{5}=0.
\end{eqnarray*}
 We chose $b_{1}=32$ and $o_{11}=6$, which allowed us to write 
\[
o_{5}+o_{7}+o_{9}=26.
\]
 In the system of equations, if we add the expressions together then
we get 
\[
-b_{2}+b_{3}-b_{4}+\left(\frac{32}{10}-\frac{40}{10}+\frac{18}{10}\right)\cdot o_{9}+\left(\frac{6}{3}-\frac{3}{3}\right)\cdot o_{7}+\frac{3}{3}\cdot o_{5}=0,
\]
 which is 
\[
-b_{2}+b_{3}-b_{4}=-\left(o_{9}+o_{7}+o_{5}\right).
\]
 From expression \ref{eq:83} it follows that 
\begin{eqnarray*}
6 & = & o_{11}=32-\left(o_{5}+o_{7}+o_{9}\right)\\
 &  & =b_{1}-b_{2}+b_{3}-b_{4}.
\end{eqnarray*}
 What this means is, if we start with the choices of $b_{1}=32$ and
$o_{11}=6$, for \textit{any} $o_{5},o_{7},o_{9}$ which satisfy $o_{5}+o_{7}+o_{9}=26$
it will follow that $b_{1}-b_{2}+b_{3}-b_{4}=6.$ For example, if
we choose 
\[
o_{5}=24,\,o_{7}=1,\,o_{9}=1,
\]
which are the \textit{wrong} values for $o_{5},o_{7},o_{9}$, the
system of equations tells us 
\begin{align*}
b_{2} & =1\cdot\frac{32}{10}=\frac{16}{5}\\
b_{3} & =1\cdot\frac{40}{10}-1\cdot2=2\\
b_{4} & =1\cdot\frac{18}{10}-1\cdot\frac{3}{3}+24=\frac{124}{5},
\end{align*}
 which yields 
\begin{equation}
O_{11}=\frac{1}{6}\cdot\left(32\cdot\left(\frac{n(n+1)}{2}\right)^{4}-\frac{16}{5}\cdot\left(\frac{n(n+1)}{2}\right)^{3}+2\cdot\left(\frac{n(n+1)}{2}\right)^{2}-\frac{124}{5}\cdot O_{5}\right).\label{eq:84}
\end{equation}
 We still have 
\[
6=32-\frac{16}{5}+2-\frac{124}{5},
\]
 but of course expressions \ref{eq:82} and \ref{eq:84} produce different
values. Therefore such an observation would have been insufficient
to suggest we had found the correct expression.
\begin{description}
\item [{Exercise}] prove the result for the general case: if $b_{1},b_{2},b_{3},\ldots,b_{m-1}$
are the coefficients for the terms in $O_{2m+1}$ and the conjectures
of (\ref{eq:83}) \textit{are} true, 
\[
b_{1}-b_{2}+b_{3}\mp\cdots\mp b_{m-1}=m+1.
\]
 For $E_{2m}$ and its coefficients, state and prove the analogous
result.\end{description}

\end{document}